\documentclass[10pt, twoside]{article}
\usepackage{fancyhdr}
\usepackage{amssymb}
\usepackage{theorem}
\usepackage{array}
\usepackage{amsmath}
\usepackage{graphicx}

\newtheorem{theorem}{Theorem}[section]
\newtheorem{proposition}[theorem]{Proposition}
\newtheorem{lemma}[theorem]{Lemma}

\newtheorem{definition}[theorem]{Definition}
\newtheorem{remark}[theorem]{Remark}

\newtheorem{observation}[theorem]{Observation}
\newtheorem{corollary}[theorem]{Corollary}

\renewcommand{\arraystretch}{2}
\newlength{\blspace} 
\newsavebox{\qedbox} 

\sbox{\qedbox}{%
\begin{picture}(5.1,5)%
\put(0,0){\framebox(5,5){}}%
\end{picture}} %
\newcommand{\qed}{\hfill$\blacksquare$}%
\newenvironment{proof}[1][]%
{\begin{trivlist}\item{\textbf{Proof#1.}\hspace{\blspace}}}%
{\qed\end{trivlist}}

\fancypagestyle{plain}{%
   \fancyhead{} 
}

\newcommand{\gln}[1]{\mathrm{GL}_{#1}\mathbb{C}}
\newcommand{\lgln}[1]{\mathfrak{gl}_{#1}\mathbb{C}}
\newcommand{\sln}[1]{\mathfrak{sl}_{#1}\mathbb{C}}

\renewcommand{\arraystretch}{1.0}

\begin{document}

\thispagestyle{plain}

\begin{center}
{\Large \textbf{A vector partition function for the multiplicities of $\sln{k}$}}
\bigskip

Sara \textsc{Billey}\footnote{Supported by NSF grant DMS-9983797.}, Victor \textsc{Guillemin} and Etienne \textsc{Rassart}\footnote{Supported by NSERC (Canada), FCAR (Qu\'ebec) and a Y. T. Li fellowship (MIT)}\\[2mm] 
\textit{Department of Mathematics, Massachusetts Institute of Technology}\\
\texttt{\{sara, vwg, rassart\}@math.mit.edu}
\bigskip

\textrm{July 15, 2003}
\end{center}
\bigskip
       
\begin{abstract}
We use Gelfand-Tsetlin diagrams to write down the weight multiplicity function for the Lie algebra $\sln{k}$ (type $A_{k-1}$) as a single partition function. This allows us to apply known results about partition functions to derive interesting properties of the weight diagrams. We relate this description to that of the Duistermaat-Heckman measure from symplectic geometry, which gives a large-scale limit way to look at multiplicity diagrams. We also provide an explanation for why the weight polynomials in the boundary regions of the weight diagrams exhibit a number of linear factors. Using symplectic geometry, we prove that the partition of the permutahedron into domains of polynomiality of the Duistermaat-Heckman function is the same as that for the weight multiplicity function, and give an elementary proof of this for $\sln{4}$ ($A_3$).  
\end{abstract}

\section{Introduction}

\newenvironment{Theorem}[1]{\noindent\textbf{Theorem #1}\itshape}{}
\newenvironment{Lemma}[1]{\noindent\textbf{Lemma #1}\itshape}{}
\newenvironment{Corollary}[1]{\noindent\textbf{Corollary #1}\itshape}{}

For a long time there has been a lot of interest, both in mathematics and physics, in finding ways to determine with what multiplicity a weight appears in the weight space decomposition of a finite-dimensional irreducible representation of a semisimple Lie algebra. Although there is a multitude of formulas to compute these multiplicities, involving partition functions (Kostant's formula), recursions (Freudenthal's formula), counting paths (Littelmann's formula), this is still a computationally hard problem. For type $A$ ($\mathrm{SL}_k\mathbb{C}$, $\mathrm{GL}_k\mathbb{C}$, $\mathrm{SU}(k)$), these multiplicities are known to be the Kostka numbers, which express the Schur symmetric functions in terms of the monomial symmetric functions.

Here we explore the structure of the weight diagrams in type $A$, not from a symmetric functions perspective, but using an array of tools from combinatorics, convex geometry and symplectic geometry, such as Gelfand-Tsetlin diagrams, Kostant's multiplicity formula, and the so-called ``Quantization Commutes with Reduction'' Theorem. We describe how the weight diagrams are partitioned into domains of polynomiality, and how this is related to the Duistermaat-Heckman function studied by symplectic geometers.

After a brief reminder about the structure of the Lie algebra $\sln{k}$, we introduce our main tools, Gelfand-Tsetlin diagrams and partition functions. Gelfand-Tsetlin theory provides a way of computing weight multiplicities by counting certain combinatorial diagrams, or equivalently, by counting the number of integer lattice points inside certain polytopes. We will use this and some notions from linear and integer programming to reduce this counting problem to evaluating a single partition function.

\begin{Theorem}{\ref{thm:SPF}}
For every $k$, we can find integer matrices $E_{k}$ and $B_{k}$ such that the multiplicity function for $\sln{k}$ can be written as 
\begin{displaymath}
m_\lambda(\beta) = \phi_{E_k}\left(B_k\left(\begin{array}{c}\lambda\\ \beta\end{array}\right)\right)\,.
\end{displaymath} 
\end{Theorem}

Expressing the multiplicities as a single partition function allows us to use general facts about partition functions and their chamber complexes to derive interesting properties of the weight diagrams. For example, the multiplicities have the following polynomiality property. 

\begin{Theorem}{\ref{thm:PolynomialityComplex}}
There is a chamber complex $\mathcal{C}^{(k)}$ on which the weight multiplicity function is determined by polynomials of degree ${k-1 \choose 2}$ in the $\beta_i$, with coefficients of degree ${k-1 \choose 2}$ in the $\lambda_j$.
\end{Theorem}

From this theorem we can deduce a pointwise scaling property (i.e. for fixed $\lambda$ and $\beta$). This property (Corollary~\ref{cor:scaling}) was known already in the context of symmetric function theory, where it was proved using a fermionic formula for the Kostka-Foulkes polynomials (see \cite{Kirillov}). It shows that although the Gelfand-Tsetlin polytopes are not always integral polytopes \cite{DLMA}, their Ehrhart quasipolynomials are in fact always polynomials.

The partition of the weight diagram into its domains of polynomiality can be described explicitly.  The convex hull of a weight diagram is a permutahedron.   There is in symplectic geometry a function on the permutahedron, called the Duistermaat-Heckman function, that approximates the weight multiplicities and is known to be piecewise polynomial. Its domains of polynomiality are convex subpolytopes of the permutahedron, and there is an explicit description of the partition in terms of walls separating the domains. Using known results on quantization and reduction of symplectic manifolds, we can prove that the Duistermaat-Heckman function and the weight multiplicity function give rise to the same partition of the permutahedron.

\begin{Theorem}{\ref{thm:SamePartition}}
The partitions of the permutahedron for $\mathfrak{su}(k)$ (or $\sln{k}$) into its domains of polynomiality for the weight multiplicities and for the Duistermaat-Heckman measure are the same.  Namely, the walls are determined by convex hulls of the form $\mathrm{conv}(W \cdot \sigma(\lambda) )$ where $\sigma \in \mathfrak{S}_k$ and $W$ is any parabolic subgroup of $\mathfrak{S}_k$ generated by all reflections corresponding to roots orthogonal to a conjugate of a fundamental weight.
\end{Theorem}

In Kostant's multiplicity formula, multiplicities are expressed as a sum of partition functions evaluated at $k!$ points shifted by a factor depending on the choice of a positive root system. We can take advantage of the apparent lack of symmetry of Kostant's multiplicity formula to find interesting factorization patterns in the weight polynomials of the boundary regions of the weight diagrams. 

\begin{Theorem}{\ref{thm:factorization}}
Let $R$ be a domain of polynomiality for the weight diagram of the irreducible representation of $\sln{k}$ with highest weight $\lambda$, and $p_R$ be its weight polynomial. Suppose that $R$ has a facet lying on the boundary of the permutahedron for $\lambda$ that has $\theta(\omega_j)$ as its normal vector, for some $\theta\in\mathfrak{S}_k$. If $\gamma = \gamma(\lambda)$ is the defining equation of the hyperplane supporting that facet,  then $p_R$ is divisible by the $j(k-j)-1$ linear factors $\gamma+1$, $\gamma+2$, $\ldots$, $\gamma+j(k-j)-1$, or $\gamma-1$, $\gamma-2$, $\ldots$, $\gamma-j(k-j)+1$.
\end{Theorem}

The main tool for proving this theorem is a family of hyperplane arrangements, called Kostant arrangements, on whose regions we have different polynomials giving the multiplicities. The Kostant arrangement also provides a method for finding linear factors in the difference between the weight polynomials of two adjacent regions. A generalization of the Kostant arrangements is also essential to the proof of Theorem~\ref{thm:PolynomialityComplex}, which establishes that although in general we get quasipolynomials in the chambers of the complex associated to a vector partition function, we get polynomials for the weight multiplicity function in type $A$.

\begin{Theorem}{\ref{thm:InternalFactorization}}
Let $R_1$ and $R_2$ be two adjacent top-dimensional domains of polynomiality of the permutahedron for a generic dominant weight $\lambda$ of $\sln{k}$, and suppose that the normal to their touching facets is in the direction $\sigma(\omega_j)$ for some $\sigma\in\mathfrak{S}_k$. If $p_1$ and $p_2$ are the weight polynomials of $R_1$ and $R_2$, and $\gamma$ is the linear functional defining the wall separating them, then the jump $p_1-p_2$ either vanishes or has the $j(k-j)-1$ linear factors 
\begin{displaymath}
(\gamma-s^{-}+1), (\gamma-s^{-}+2), \ldots, \gamma, \ldots, (\gamma+s^{+}-2), (\gamma+s^{+}-1)
\end{displaymath}
for some integers $s^{-}, s^{+} \geq 0$ satisfying
\begin{displaymath}
s^{-} + s^{+} = j(k-j)\,.
\end{displaymath}
\end{Theorem}

Similar factorization phenomena were recently observed to hold for general vector partition functions by Szenes and Vergne \cite{SzenesVergne}.

Finally, we explicitly compute the chamber complex for $A_3$, and find it is not optimal, but that we can glue together parts of it to obtain a simpler complex. We can deduce symbolically from the form of this complex that the optimal partitions of the permutahedron for $A_3$ under the weights and the Duistermaat-Heckman measure are the same. Computing the chamber complex for $A_{3}$ is nontrivial because of the complexity of the arrangement.  To the best of our knowledge, these computations for generic dominant weights of $A_3$ have not been done. A study was done by Guillemin, Lerman and Sternberg in \cite{GLS} for some of the degenerate cases when $\lambda$ has a nontrivial stabilizer. The number of domains of polynomiality turns out to be significantly larger than they originally suspected.

\subsection{The Lie algebra $\sln{k}$ (type $A_{k-1}$)}

The simple Lie algebra $\sln{k}$ is the subalgebra of $\lgln{k} \cong \mathrm{End}(\mathbb{C}^k)$ consisting of traceless $k\times k$ matrices over $\mathbb{C}$. We will take as its Cartan subalgebra $\mathfrak{h}$ its subspace of traceless diagonal matrices. The roots and weights live in the dual $\mathfrak{h}^*$ of $\mathfrak{h}$, which can be identified with the subspace $x_1+\cdots+x_{k}=0$ of $\mathbb{R}^{k}$. The roots are $\{e_i-e_j\ : \ 1\leq i\neq j\leq k\}$, and we will choose the positive ones to be $\Delta_+ = \{e_i-e_j\ : \ 1\leq i<j\leq k\}$. The simple roots are then $\alpha_i = e_i-e_{i+1}$, for $1\leq i\leq k-1$, and for these simple roots, the fundamental weights are
\begin{equation}
\omega_i = \frac{1}{k}(\underbrace{k-i,k-i,\ldots,k-i}_{\textrm{$i$ times}},\underbrace{-i,-i,\ldots,-i}_{\textrm{$k-i$ times}})\,, \qquad 1\leq i\leq k-1\,.
\end{equation}
The fundamental weights are defined such that $\langle\alpha_i,\omega_j\rangle=\delta_{ij}$, where $\langle\cdot,\cdot\rangle$ is the usual dot product. The integral span of the simple roots and the fundamental weights are the root lattice $\Lambda_R$ and the weight lattice $\Lambda_W$ respectively. The root lattice is a finite index sublattice of the weight lattice, with index $k-1$. 

For our choice of positive roots, 
\begin{equation}
\delta = \frac{1}{2}\sum_{\alpha\in\Delta_+}\alpha = \sum_{j=1}^{k-1}\omega_k = \frac{1}{2}(k-1,k-3,\ldots,-(k-3),-(k-1))\,.
\end{equation}

The Weyl group for $\sln{k}$ is the symmetric group $\mathfrak{S}_k$ acting on $\{e_1,\ldots,e_k\}$ (i.e. $\sigma(e_i)=e_{\sigma(i)}$), and with the choice of positive roots we made, the fundamental Weyl chamber will be $C_0 = \{(\lambda_1,\ldots,\lambda_k)\ : \ \sum_{i=1}^k\lambda_i=0 \ \textrm{and } \lambda_1\geq\cdots\geq\lambda_k\}$. The action of the Weyl group preserves the root and weight lattices. The \emph{Weyl orbit} of a weight $\lambda$ is the set $\mathfrak{S}_k\cdot\lambda=\{\sigma(\lambda)\ :\ \sigma\in\mathfrak{S}_k\}$. We refer to the convex hull of $\mathfrak{S}_k\cdot\lambda$ as the \textit{permutahedron} associated to $\lambda$. Weights lying in the fundamental Weyl chamber are called \emph{dominant}, and we will call elements of the Weyl orbits of the fundamentals weights \emph{conjugates of fundamental weights}.

The finite dimensional representations of $\sln{k}$ are indexed by the dominant weights $\Lambda_W\cap C_0$, and for a given dominant weight $\lambda$, there is a unique irreducible representation $\rho_\lambda : \sln{k} \rightarrow \mathfrak{gl}(V_\lambda)$ with highest weight $\lambda$, up to isomorphism. Details about their construction are well-known and can be found in \cite{Fulton} or \cite{FultonHarris}, for example. We have the weight space decomposition according to the action of $\mathfrak{h}$
\begin{equation}
V_\lambda = \bigoplus_{\beta}\big(V_\lambda\big)_\beta\,.
\end{equation}

The weights of this representation (those $\beta$'s for which $\big(V_\lambda\big)_\beta\neq0$) are finite in number, and they can be characterized as follows (see \cite{Humphreys}): they are exactly the points $\beta$ of the weight lattice $\Lambda_W$ that lie within the convex hull of the orbit of $\lambda$ under the Weyl group action, denoted $\mathrm{conv}(\mathfrak{S}_k\cdot\lambda)$, and such that $\lambda-\beta$ lies in the root lattice. Hence
\begin{equation}
V_\lambda = \bigoplus_{\beta\in (\lambda+\Lambda_R)\cap\mathrm{conv}(\mathfrak{S}_k\cdot\lambda)}\big(V_\lambda\big)_\beta\,.
\end{equation}

The multiplicity $m_\lambda(\beta)$ of the weight $\beta$ in $V_\lambda$ is the dimension of $(V_\lambda\big)_\beta$, and all the conjugates of $\beta$ under $\mathfrak{S}_k$ have the same multiplicity. The \textit{weight diagram} of $V_\lambda$ consists of the weights of $V_\lambda$ (as a subset of $\Lambda_W$) together with the data of their multiplicities.

There are several ways to compute weight multiplicities. An important one is Kostant's multiplicity formula \cite{Kostant}, which can be deduced from Weyl's character formula (see \cite{Humphreys,Stembridge1}). We first need to define the Kostant partition function given a choice of positive root system $\Delta_+$:
\begin{equation}
K(v) = \Big|\Big\{(k_\alpha)_{\alpha\in\Delta_+}\in\mathbb{N}^{|\Delta_+|}\ : \ \sum_{\alpha\in\Delta_+}k_\alpha\alpha = v\Big\}\Big|\,,
\end{equation}
i.e. $K(v)$ is the number of ways that $v\in\mathfrak{h}^*$ can be written as a sum of positive roots.

Kostant's multiplicity formula \cite{Kostant} is then
\begin{equation}
\label{eqn:KMF}
m_\lambda(\beta) = \sum_{\sigma\in\mathfrak{S}_k}(-1)^{l(\sigma)}K(\sigma(\lambda+\delta)-(\beta+\delta))\,,
\end{equation}
where $l(\sigma)$ is the number of inversions $\sigma$. Kostant's partition function and multiplicity formula extend to all complex semisimple Lie algebras.  See \cite{Humphreys} for more details.

\subsection{Gelfand-Tsetlin diagrams}

Gelfand-Tsetlin diagrams were introduced by Gelfand and Tsetlin \cite{GelfandTsetlin} as a way to index the one-dimensional subspaces of the (polynomial) representations of $\gln{k}$. Their construction relies on a theorem of Weyl that describes how the restriction to $\gln{k-1}$ of an irreducible representation of $\gln{k}$ breaks down into irreducible representations of $\gln{k-1}$ (see \cite{BarceloRam,GelfandTsetlin,Zelobenko}). They are equivalent to semistandard tableaux (see \cite{GoodmanWallach}), but they have a ``linear'' structure that we will exploit. 

\begin{definition}
Let $\nu=(\nu_1,\ldots,\nu_m)$ and $\gamma=(\gamma_1,\ldots,\gamma_{m-1})$ be two partitions. We will say that $\gamma$ \emph{interlaces} $\nu$, and write $\gamma \lhd \nu$, if
\begin{displaymath}
\nu_1 \geq \gamma_1 \geq \nu_2 \geq \gamma_2 \geq \nu_3 \geq \cdots \geq \nu_{m-1} \geq \gamma_{m-1} \geq \nu_m\,.
\end{displaymath}
\end{definition}

\begin{theorem}{\rm\textbf{(Weyl's branching rule \cite{GoodmanWallach,Zelobenko})}}
Let $\rho_\lambda$ be the (polynomial) irreducible representation of $\gln{k}$ with highest weight $\lambda = \lambda_1\geq\lambda_2\geq\ldots\lambda_k\geq0$. The decomposition of the restriction of $\rho_\lambda$ to $\gln{k-1}$ into irreducible representations of $\gln{k-1}$ is given by 
\begin{equation}
{\rho_\lambda}_{\big|_{\gln{k-1}}} = \bigoplus_{\mu \lhd \lambda}\rho_\mu\,.
\end{equation} 
\end{theorem}

After restricting $\rho_\lambda$ to $\gln{k-1}$ and breaking it into $\gln{k-1}$-irreducibles, we can restrict to $\gln{k-2}$:
\begin{equation}
{\rho_\lambda}_{\big|_{\gln{k-2}}} = \Big({\rho_\lambda}_{\big|_{\gln{k-1}}}\Big)_{\big|_{\gln{k-2}}} = \Big(\bigoplus_{\mu \lhd \lambda}\rho_\mu\Big)_{\big|_{\gln{k-2}}} = \bigoplus_{\mu \lhd \lambda}\Big({\rho_\mu}_{\big|_{\gln{k-2}}}\Big)\,.
\end{equation}

Again, we can apply Weyl's branching rule to each $\rho_\mu$ to break them into irreducible representations of $\gln{k-2}$ to get
\begin{equation}
{\rho_\lambda}_{\big|_{\gln{k-2}}} = \bigoplus_{\nu\lhd\mu\lhd\lambda}\rho_\nu\,.
\end{equation}

We can keep going recursively, and for convenience, let us denote by $\lambda^{(m)}=\lambda^{(m)}_1\geq\cdots\geq\lambda^{(m)}_m\geq0$ the partitions indexing the irreducible representations of $\gln{m}$. We then get that
\begin{equation}
{\rho_\lambda}_{\big|_{\gln{1}}} = \bigoplus_{\lambda^{(1)}\lhd\cdots\lhd\lambda^{(k)}=\lambda}\rho_{\lambda^{(1)}}\,.
\label{eqn:GTdecomp}
\end{equation} 

\newcommand{\dddots}{\raisebox{-1.1mm}[0cm][0cm]{$\cdot$}\,\raisebox{0mm}[0cm][0cm]{$\cdot$}\,\raisebox{1.1mm}[0cm][0cm]{$\cdot$}}

\begin{definition}
A sequence of partitions of the form ${\lambda^{(1)}\,\lhd\,\cdots\,\lhd\,\lambda^{(k)}\,=\,\lambda}$ is called a \emph{Gelfand-Tsetlin diagram} for $\lambda$, and can be viewed schematically as
\begin{equation}
\renewcommand{\arraystretch}{1.4}
\begin{array}{c@{\hspace{2mm}}c@{\hspace{2mm}}c@{\hspace{2mm}}c@{\hspace{2mm}}c@{\hspace{2mm}}c@{\hspace{2mm}}c@{\hspace{2mm}}c@{\hspace{2mm}}c}
\lambda^{(k)}_1 & & \lambda^{(k)}_2 & & \cdots & & \lambda^{(k)}_{k-1} & & \lambda^{(k)}_k \\
& \lambda^{(k-1)}_{1} & & \lambda^{(k-1)}_{2} & & \cdots & & \lambda^{(k-1)}_{k-1} & \\
& & \ddots & & \vdots & & \ \raisebox{0.5mm}{\dddots}  & & \\ 
& & & \lambda^{(2)}_{1} & & \lambda^{(2)}_{2} & & & \\
& & & & \lambda^{(1)}_1 & & & & 
\end{array}
\end{equation}
with $\lambda^{(k)}_j=\lambda_j$ and each $\lambda^{(i)}_j$ is a nonnegative integer satisfying
\begin{equation}
\lambda^{(i+1)}_{j}\  \geq \ \lambda^{(i)}_{j} \ \geq \ \lambda^{(i+1)}_{j+1} \qquad\textrm{for \,$1\leq i\leq k-1$, $1\leq j\leq i$\,.}
\label{eqn:GTineq}
\end{equation}
\end{definition}
\smallskip

Let $V_\mathcal{D}$ be the one-dimensional subspace of $V_{\lambda}$ corresponding to a Gelfand-Tsetlin diagram $\mathcal{D}$. It is shown in \cite{Zelobenko} that $V_\mathcal{D}$ lies completely within one weight space in the weight space decomposition of $V_\lambda$:  $V_\mathcal{D}\subseteq\big(V_\lambda\big)_{\beta}$ if 
\begin{equation}
\beta_m = \sum_{i=1}^{m}\lambda^{(m)}_i-\sum_{i=1}^{m-1}\lambda^{(m-1)}_i \qquad \textrm{for $1\leq m\leq k$}
\label{eqn:GTweight1}
\end{equation}
or, equivalently,
\begin{equation}
\beta_1 + \cdots + \beta_m = \sum_{i=1}^{m}\lambda^{(m)}_i \qquad \textrm{for $1\leq m\leq k$}.
\label{eqn:GTweight2}
\end{equation}

Hence Gelfand-Tsetlin diagrams for $\lambda$ correspond to the same weight if all their row sums are the same. This discussion is summarized in the following theorem due to Gelfand, Tsetlin and Zelobenko.

\begin{theorem}{\rm\textbf{(\cite{GelfandTsetlin,Zelobenko})}}
For $\lambda=(\lambda_1,\ldots,\lambda_k)$, the number of Gelfand-Tsetlin diagrams with first row $\lambda$ is the dimension of the irreducible representation $V_\lambda$ of $\gln{k}$ with highest weight $\lambda$. Furthermore, the multiplicity $m_\lambda(\beta)$ of the weight $\beta$ in the irreducible representation of $\gln{k}$ with highest weight $\lambda$ is given by the number of Gelfand-Tsetlin diagrams with first row $\lambda$ such that equation~(\ref{eqn:GTweight1}) (or (\ref{eqn:GTweight2})) is satisfied. 
\end{theorem}

Two irreducible representations $V_\lambda$ and $V_{\gamma}$ of $\lgln{k}$ restrict to the same irreducible representation of $\sln{k}$ if $\lambda_i-\gamma_i$ is some constant independent of $i$ for all $i$. Hence we might as well require that the $\lambda_i$ sum up to zero. However, normalizing the sum this way can introduce fractional values of $\lambda$, so we'll have to translate $\lambda$ back to integer values when writing down Gelfand-Tsetlin diagrams for those representations, or, equivalently, translate the integer lattice along with $\lambda$, so that the inequalities
\begin{displaymath}
\lambda^{(i+1)}_{j}\  \geq \ \lambda^{(i)}_{j} \ \geq \ \lambda^{(i+1)}_{j+1} \qquad\textrm{for \,$1\leq i\leq k-1$, $1\leq j\leq i$\,,}
\end{displaymath}
always have
\begin{displaymath}
\lambda^{(i+1)}_{j} - \lambda^{(i)}_{j} \in \mathbb{N} \qquad \textrm{and} \qquad \lambda^{(i)}_{j} - \lambda^{(i+1)}_{j+1} \in \mathbb{N}\,.
\end{displaymath}

There is a geometrical way to view the enumeration of the number of Gelfand-Tsetlin diagrams for a given $\lambda$. With $\lambda^{(k)}=\lambda$ fixed, we can let all the other variables $\{\lambda^{(m)}_i\ : \ 1\leq i\leq m, 1\leq m<k\}$ be real variables. The system of inequalities~(\ref{eqn:GTineq}) among the entries of Gelfand-Tsetlin diagrams, when viewed over the reals, defines a rational polytope, called the \emph{Gelfand-Tsetlin polytope for $\lambda$} and denoted $\mathrm{GT}_\lambda$. $\mathrm{GT}_\lambda$ has dimension at most ${k\choose 2}$, and equal to that number if the $\lambda_i$'s are distinct. We can consider the intersection of this polytope with the affine subspace obtained by fixing a weight $\beta$ (fixing the row sums using equations~(\ref{eqn:GTweight1}) or~(\ref{eqn:GTweight2})). We also get a rational polytope this way, called the \emph{Gelfand-Tsetlin polytope for $\lambda$ and $\beta$} and denoted $\mathrm{GT}_{\lambda,\beta}$\,. Its dimension is at most ${k-1\choose 2}$. Kirillov conjectured in \cite{Kirillov} that the polytopes $\mathrm{GT}_{\lambda,\beta}$ are integral polytopes, but this was recently disproved by De~Loera and McAllister \cite{DLMA}.

The upshot is that integer solutions to the Gelfand-Tsetlin diagram constraints then translate into integer points inside the polytopes, hence the number of Gelfand-Tsetlin diagrams of weight $\beta$ for $\lambda$ is the number of integer points in the polytope $\mathrm{GT}_{\lambda,\beta}$.

\subsection{Partition functions and chamber complexes}

Partition functions arise in the representation theory of the semisimple Lie algebras through Kostant's formula for the multiplicities \eqref{eqn:KMF}.   Kostant's partition function sends a vector in the root lattice to the number of ways it can be written down as a linear combination with nonnegative integer coefficients of the positive roots, and this is a simple example of a more general class of functions, called \emph{vector partition functions}.

\begin{definition}
Let $M$ be a $d \times n$ matrix over the integers, such that $\mathrm{ker}M \cap \mathbb{R}^n_{\geq 0} = 0$. The \emph{vector partition function} (or simply \emph{partition function}) associated to $M$ is the function
\begin{displaymath}
\begin{array}{rccl}
\phi_M : & \mathbb{Z}^d & \longrightarrow & \mathbb{N}\\
& b & \mapsto & |\{x \in \mathbb{N}^n \ : \  Mx = b\}|
\end{array}
\end{displaymath}
\end{definition}

The condition $\mathrm{ker}M \cap \mathbb{R}^n_{\geq 0} = 0$ forces the set $\{x \in \mathbb{N}^n \ : \  Mx = b\}$ to have finite size, or equivalently, the set $\{x \in \mathbb{R}^n_{\geq 0} \ : \  Mx = b\} $ to be compact, in which case it is a polytope $P_b$, and the partition function is the number of integral points (lattice points) inside it.

Also, if we let $M_1, \ldots, M_n$ denote the columns of $M$ (as column-vectors), and $x=(x_1, \ldots, x_n) \in \mathbb{R}^n_{\geq 0}$, then $Mx = x_1M_1 + x_2M_2 + \cdots + x_nM_n$ and for this to be equal to $b$, $b$ has to lie in the cone $\mathrm{pos}(M)$ spanned by the vectors $M_i$. So $\phi_M$ vanishes outside of $\mathrm{pos}(M)$. 

It is well-known that partition functions are piecewise quasipolynomial, and that the domains of quasipolynomiality form a complex of convex polyhedral cones, called the \emph{chamber complex}. Sturmfels gives a very clear explanation in \cite{Sturmfels} of this phenomenon. The explicit description of the chamber complex is due to Alekseevskaya, Gel'fand and Zelevinski$\breve{\i}$ \cite{AGZ}. There is a special class of matrices for which partition functions take a much simpler form. Call an integer $d\times n$ matrix $M$ of full rank $d$ \emph{unimodular} if every nonsingular $d\times d$ submatrix has determinant $\pm 1$. For unimodular matrices, the chamber complex determines domains of polynomiality instead of quasipolynomiality \cite{Sturmfels}.

It is useful for what follows to describe how to obtain the chamber complex of a partition function. Let $M$ be a $d\times n$ integer matrix of full rank $d$ and $\phi_M$ its associated partition function. For any subset $\sigma \subseteq \{1,\ldots,n\}$, denote by $M_\sigma$ the submatrix of $M$ with column set $\sigma$, and let $\tau_\sigma = \mathrm{pos}(M_\sigma)$, the cone spanned by the columns of $M_\sigma$. Define the set $\mathcal{B}$ of \emph{bases} of $M$ to be
\begin{displaymath}
\mathcal{B} = \{\sigma \subseteq \{1,\ldots,n\}\ : \ |\sigma| = d \ \ \textrm{and } \ \mathrm{rank}(M_\sigma) = d\}\,.
\end{displaymath}
$\mathcal{B}$ indexes the invertible $d\times d$ submatrices of $M$. The \emph{chamber complex} of $\phi_M$ is the common refinement of all the cones $\tau_\sigma$, as $\sigma$ ranges over $\mathcal{B}$ (see \cite{AGZ}). A theorem of Sturmfels \cite{Sturmfels} describes exactly how partition functions are quasipolynomial over the chambers of that complex.

\subsection{The chamber complex for the Kostant partition function}
\label{subsec:KPF}

If we let $M_{A_n}$ be the matrix whose columns are the positive roots $\Delta_+^{(A_n)}$ of $A_n$, written in the basis of simple roots, then we can write Kostant's partition function in the matrix form defined above as
\begin{displaymath}
K_{A_n}(v) = \phi_{M_{A_n}}(v)\,.
\end{displaymath}

Consider for example the simple Lie algebra $\sln{4}$, or $A_3$. The positive roots are $\Delta^{(A_3)}_+ = \{e_i-e_j\ :\ 1\leq i<j\leq 4\}$. Writing the positive roots in the basis of simple roots, we have $\Delta_+^{(A_3)} = \{\alpha_1,\alpha_2,\alpha_3,\alpha_1+\alpha_2, \alpha_2+\alpha_3, \alpha_1+\alpha_2+\alpha_3\}$. This gives
\begin{displaymath}
\renewcommand{\arraystretch}{1.0}
M_{A_3} = \left(\begin{array}{c@{\hspace{4mm}}c@{\hspace{4mm}}c@{\hspace{4mm}}c@{\hspace{4mm}}c@{\hspace{4mm}}c}
1  &  0 &  0  &  1 & 0  & 1 \\ 
0  &  1 &  0  &  1 & 1  & 1 \\
0  &  0 &  1  &  0 & 1 &  1
\end{array}\right)
\renewcommand{\arraystretch}{1.0}
\end{displaymath}
which has the bases
\begin{displaymath}
\mathcal{B} = \{123, 125, 126, 134, 135, 136, 145, 146, 234, 236, 245, 246, 256, 345, 356, 456\}\,,
\end{displaymath}
where we're writing $i_1i_2i_3$ for $\{i_1, i_2, i_3\}$.

\begin{figure}[!ht]
\begin{center}
\includegraphics[width=0.8\textwidth]{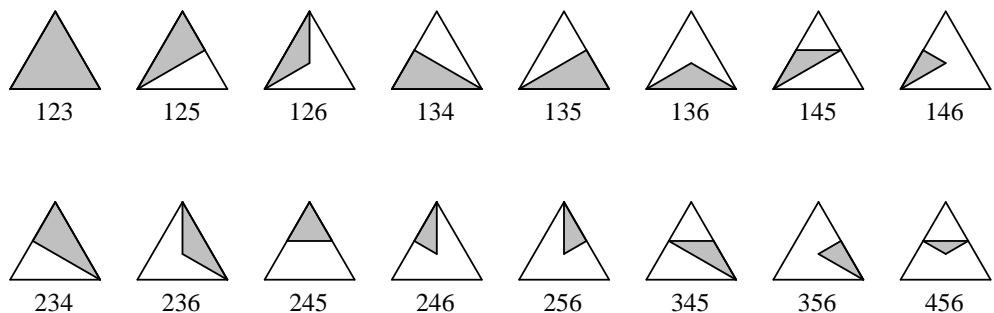}
\caption{Basis cones for the Kostant partition function of $A_3$.}
\label{fig:kostant_basis_cones}
\end{center}
\end{figure}

\begin{figure}[!ht]
\begin{center}
\includegraphics[width=0.3\textwidth]{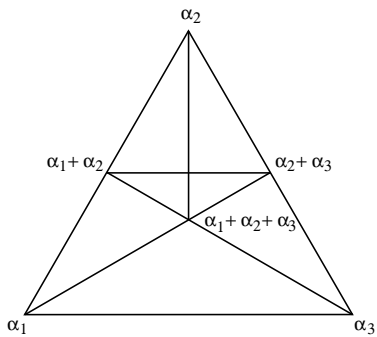}
\caption{Chamber complex for the Kostant partition function of $A_3$.}
\label{fig:kostant_complex}
\end{center}
\end{figure}

All the cones corresponding to these bases are contained in the first cone with basis $\{1,2,3\}$ which is just the positive octant in $\mathbb{R}^3$. To picture the chamber complex, we can look at the intersection of these cones with the hyperplane $x+y+z=1$. Figure~\ref{fig:kostant_basis_cones} shows the cones given by the bases of $\mathcal{B}$, while Figure~\ref{fig:kostant_complex} shows their common refinement (this originally appeared in \cite{DeLoeraSturmfels}). Finally, since it is readily checked that $M_{A_3}$ is unimodular, this shows that the Kostant partition function for $A_3$ has 7 domains of polynomiality. It is an open problem mentioned by Kirillov in \cite{Kirillov} to determine the numbers of chambers for the Kostant partition functions for the Lie algebras $A_n$. De~Loera and Sturmfels \cite{DeLoeraSturmfels} have computed the numbers for $n \leq 6$ and computed the polynomial associated to each chamber for $n \leq 5$.

The following lemma is a well-known fact about $M_{A_n}$ and can be deduced from general results on matrices with columns of $0$'s and $1$'s where the $1$'s come in a consecutive block (see \cite{Schrijver}).

\begin{lemma}
The matrix $M_{A_n}$ is unimodular for all $n$. 
\end{lemma}

$M_{A_n}$ unimodular means that the Kostant partition functions for $A_n$ is polynomial instead of quasipolynomial on the cells of the chamber complex. In general, for $M$ unimodular, the polynomial pieces have degree at most the number of columns of the matrix minus its rank (see \cite{Sturmfels}). In our case, $M_{A_n}$ has rank $n$ and as many columns as $A_n$ has positive roots, ${n+1\choose 2}$. Hence the Kostant partition function for $A_n$ is piecewise polynomial of degree at most ${n+1\choose 2}-n = {n\choose 2}$.

\begin{remark}
\label{rem:LocallyPolynomial}
In view of Kostant's formula for the weight multiplicities \eqref{eqn:KMF}, this means that the multiplicity function $m_\lambda(\beta)$ for $A_n$ is piecewise polynomial of degree at most ${n\choose 2}$ in the $\beta$-coordinates if the $\lambda$-coordinates are fixed, or degree ${n\choose 2}$ in the $\lambda$-coordinates if the $\beta$-coordinates are fixed. So we can regard it as a piecewise polynomial function of degree ${n\choose 2}$ in the $\beta_i$'s, with coefficients of degree ${n\choose 2}$ in the $\lambda_j$'s. This will be made precise in Sections~\ref{sec:KostantArrangements} and~\ref{sec:PolynomialityComplex}
\end{remark}

From now on, we will be interested in the multiplicity function for $\sln{k}$, of type $A_{k-1}$, and thus use the results above with $n=k-1$.

\section{The multiplicity function as a single partition function}
\label{sec:SPF}

Our first theorem presents a new conceptual approach to computing multiplicities. This approach is efficient for large $\lambda$ in low ranks. It has the additional advantages of allowing us to use known facts about partitions functions. 

\begin{theorem}
\label{thm:SPF}
For every $k$, we can find integer matrices $E_{k}$ and $B_{k}$ such that the multiplicity function for $\sln{k}$ can be written as 
\begin{equation}
m_\lambda(\beta) = \phi_{E_k}\left(B_k\left(\begin{array}{c}\lambda\\ \beta\end{array}\right)\right)\,.
\end{equation} 
\end{theorem}

\begin{proof}
Consider a Gelfand-Tsetlin diagram where we will think of $\lambda=(\lambda_1,\ldots,\lambda_k)$ and $\beta=(\beta_1,\ldots,\beta_k)$ as parameters, with the conditions that $\sum_{i=1}^k\lambda_i = \sum_{i=1}^k\beta_i = 0$. The variables in the diagram are $\lambda^{(i)}_j$ with $1\leq i\leq k-1$, $1\leq j\leq i$. Each of these ${k \choose 2}$ variables is wedged between two entries at the level above, so we get a system of $2{k \choose 2} = k(k-1)$ inequalities. Using equation~(\ref{eqn:GTweight2}), relating the row sums to the $\beta_i$'s, we can get rid of the $k-1$ variables $\lambda^{(1)}_{1}$, $\lambda^{(2)}_{2}$, ..., $\lambda^{(k-1)}_{k-1}$. 
\begin{equation}
\renewcommand{\arraystretch}{1.4}
\begin{array}{c@{\hspace{2mm}}c@{\hspace{2mm}}c@{\hspace{2mm}}c@{\hspace{2mm}}c@{\hspace{2mm}}c@{\hspace{2mm}}c@{\hspace{2mm}}c@{\hspace{2mm}}c@{\hspace{2mm}}c@{\hspace{2mm}}cl}
\lambda_1 & & \lambda_2 & & \cdots & & \cdots & & \lambda_{k-1} & & \lambda_k & \quad (\beta_1 + \cdots + \beta_k = 0)\\
& \fbox{$\lambda^{(k-1)}_{1}$} & & \fbox{$\lambda^{(k-1)}_{2}$} & & \cdots & & \fbox{$\lambda^{(k-1)}_{k-2}$} & & \lambda^{(k-1)}_{k-1} & & \quad (\beta_1 + \cdots + \beta_{k-1})\\
& & \ddots & & \vdots & & \vdots & & \ \raisebox{0.5mm}{\dddots}  & & & \\ 
& & & \fbox{$\lambda^{(3)}_{1}$} & & \fbox{$\lambda^{(3)}_{2}$} & & \lambda^{(3)}_{3} & & & & \quad (\beta_1+\beta_2+\beta_3)\\
& & & & \fbox{$\lambda^{(2)}_{1}$} & & \lambda^{(2)}_{2} & & & & & \quad (\beta_1+\beta_2)\\
& & & & & \lambda^{(1)}_1 & & & & & & \quad (\beta_1) 
\end{array}
\end{equation}

The remaining variables (boxed in the above diagram) are $\lambda^{(i)}_j$ with $1\leq i\leq k-1$, $1\leq j\leq i-1$ and there are ${k-1 \choose 2}$ of them. To get a system in partition function form, we need to transform the inequalities into equalities satisfied by nonnegative variables, however the $\lambda^{(i)}_j$ can take negative values. Let
\begin{displaymath}
s^{(i)}_j = \lambda^{(i)}_j - \lambda^{(i+1)}_{j+1} \qquad 1\leq i\leq k-1, \quad 1\leq j\leq i-1
\end{displaymath}
be the differences between the variables and the ones immediately above and to the right of them, recalling that $\lambda^{(k)}_j = \lambda_j$. Upon doing the substitution in the system of inequalities, ${k-1 \choose 2}$ of the inequalities simply become $s^{(i)}_j\geq 0$ because of equation~(\ref{eqn:GTineq}). So we are left with a system of $N = k(k-1)-{k-1 \choose 2} = \frac{1}{2}(k-1)(k+2)$ inequalities in the $K = {k-1 \choose 2}$ \emph{nonnegative and integral} variables $s^{(i)}_j$, which we will relabel $s_1,\ldots,s_{K}$ for convenience. 

The final step is to transform the inequalities into equalities. To this effect, we write each inequality in the form
\begin{displaymath}
a_{m1}s_1 + a_{m2}s_2 + \cdots + a_{mK}s_K \ \leq \ \sum_{j=1}^kb_{mj}\lambda_j + \sum_{j=1}^kc_{mj}\beta_j\,,
\end{displaymath}
for $1\leq m\leq N$ and integers $a_{m1},\ldots,a_{mK}$, $b_{mj}$, $c_{mj}$ ($1\leq j\leq k$).

We introduce a slack variable for each inequality to turn it into an equality:
\begin{displaymath}
a_{m1}s_1 + a_{m2}s_2 + \cdots + a_{mK}s_K + s_{K+m} \ = \ \sum_{j=1}^kb_{mj}\lambda_j + \sum_{j=1}^kc_{mj}\beta_j\,.
\end{displaymath}

The slack variables $s_{K+1},\ldots,s_{K+N}$ are nonnegative, just like the previous $K$ $s_i$, and integral solutions to the system of inequalities will correspond to integral solutions to this system of equalities, so $s_{K+1},\ldots,s_{K+N}$ are not only nonnegative but integral.

Finally, we can write the system of equalities in matrix form:
\begin{equation}
\underbrace{\begin{array}{@{}c@{}}\left(\begin{array}{ccc|ccc}
a_{11} & \cdots & a_{1K} &  &  & \\
\vdots & \ddots & \vdots &  & I_{N} & \\
a_{N1} & \cdots & a_{NK} &  &  &
\end{array}\right)\\ \phantom{a}\end{array}}_{\displaystyle E_k}\left(\begin{array}{c}s_1\\ \vdots\\ s_K\\ s_{K+1}\\ \vdots\\ s_{K+N}\end{array}\right) = \underbrace{\begin{array}{c}\left(\begin{array}{c}\sum_{j=1}^kb_{1j}\lambda_j + \sum_{j=1}^kc_{1j}\beta_j\\[2mm] \vdots\\[2mm] \sum_{j=1}^kb_{Nj}\lambda_j + \sum_{j=1}^kc_{Nj}\beta_j\end{array}\right)\\ \phantom{a}\end{array}}_{\displaystyle B_k\left(\begin{array}{c}\lambda\\ \beta\end{array}\right)}
\end{equation}

The result follows, since the number $m_\lambda(\beta)$ of integral solutions to the Gelfand-Tsetlin inequalities is the number of all integral nonnegative solutions to this matrix system.
\end{proof}

The partition function $\phi_{E_{k}}$ in the above theorem lives on a larger dimensional space than the one we need. It takes values in $\mathbb{R}^{N} = \mathbb{R}^{(k-1)(k+2)/2}$, whereas the part that interests us, the space given by $B_k{\lambda\choose\beta}$ as the $\lambda_i$ and $\beta_j$ range over $\mathbb{R}$, has dimension $2k-2$. Let 
\begin{equation}
\tilde{B} = \left\{B_k{\lambda\choose\beta}\ : \ \lambda\in\mathbb{R}^k, \beta\in\mathbb{R}^k, \sum_{i=1}^k\lambda_i=\sum_{i=1}^k\beta_i=0\right\}\,,
\end{equation} 
then the only part of the chamber complex that is relevant to the multiplicity function is its intersection with $\tilde{B}$. Since the chamber complex is obtained as the common refinement of the base cones, we will get the same thing if we find the refinement of the base cones and then intersect the result with $\tilde{B}$, or intersect the base cones with $\tilde{B}$ first and then find the common refinement of those restricted base cones. Since we only need the restricted chamber complex, this simplifies the computation because we have to deal with $2k-2$-dimensional cones instead of $(k-1)(k+2)/2$-dimensional ones. Another bonus we get from working on $\tilde{B}$ is that on this space, $B_k$ is an invertible transformation, so we can rectify the cones to $(\lambda,\beta)$-coordinates. In effect, we remove the coordinate ``twist'' due to matrix $B_{k}$. 

\begin{definition}
We will denote by $\mathcal{C}^{(k)}$ this rectified $(2k-2)$-dimensional complex in $(\lambda,\beta)$-coordinates. 
\end{definition}

Because $E_k$ is not unimodular in general, the associated partition function will be quasipolynomial on the cells of the chamber complex. We will prove in Section~\ref{sec:PolynomialityComplex} that it is actually polynomial on the cells of the complex. As such, we will from now on refer to the domains of quasipolynomiality of the multiplicity function as domains of polynomiality.

\begin{remark}
The multiplicity function also satisfies another sort of polynomiality property. There are many ways to think of fixed type $A$ dominant weights $\lambda$ and $\beta$ as living in $\sln{r}$ for any sufficiently large $r$. It is known (see for example \cite{BKLS,KingPlunkett}) that if $m^{(r)}_\lambda(\beta)$ is the multiplicity of $\beta$ in the irreducible representation $V_\lambda$ of $\sln{r}$, then $m^{(r)}_\lambda(\beta)$ is given by a polynomial function in $r$, for $r$ large enough. Bounds on the degree of this polynomial are also given. This result is shown to extend to the other classical groups \cite{KingPlunkett} and also the classical affine Kac-Moody algebras \cite{BKLS}. In our investigation of the weight multiplicities, we instead fix the rank of the Lie algebra and study the polynomial dependence in the $\lambda$ and $\beta$ variables.
\end{remark}

\begin{definition}
For every $\lambda$ in the fundamental Weyl chamber, let
\begin{equation}
\label{eqn:Llambda}
L(\lambda) = \{(\lambda_1,\ldots,\lambda_k,\beta_1,\ldots,\beta_k)\ : \ \beta_i\in\mathbb{R}\}\,.  
\end{equation}
Note that this space is really $(k-1)$-dimensional since $\sum_j\beta_j=0$. Define also the projection
\begin{equation}
p_{\Lambda}\,:\ (\lambda_1,\ldots,\lambda_k,\beta_1,\ldots,\beta_k) \longmapsto (\lambda_1,\ldots,\lambda_k)\,.
\end{equation}
\end{definition}

\begin{remark}
\label{rem:ComplexPartition}
The intersection of $\mathcal{C}^{(k)}$ with $L(\lambda)$ will give domains of polynomiality for the weight diagram of the irreducible representation of $\sln{k}$ with highest weight $\lambda$. The partition into domains that we get this way, however, is not optimal, as shown for $\sln{4}$ in Section~\ref{sec:A2A3}. Some adjacent regions have the same weight polynomial and their union is again a convex polytope, so they can be glued together to yield a larger domain. 
\end{remark}

\begin{corollary}
\label{cor:LambdaComplex}
Let $\mathcal{C}^{(k)}_\Lambda$ be the chamber complex given by the common refinement of the projections $p_{\Lambda}(\tau)$ of the cones of $\mathcal{C}^{(k)}$ onto $\mathbb{R}^{k}$. Then $\mathcal{C}^{(k)}_\Lambda$ classifies the $\lambda$'s, in the sense that if $\lambda$ and $\lambda'$ belong to the same cell of $\mathcal{C}^{(k)}_\Lambda$, then all their domains are indexed by the same subsets of cones from $\mathcal{C}^{(k)}$, and therefore have the same corresponding polynomials. 
\end{corollary}

\begin{proof}
We can index the top-dimensional domains by the top-dimensional cones $\tau$ of $\mathcal{C}^{(k)}$. The domain indexed by cone $\tau$ is present in the weight diagram (permutahedron) for $\lambda$ if and only if $\lambda\in p_{\Lambda}(\tau)$. 
\end{proof}

\section{Domains of polynomiality via Duistermaat-Heckman theory}
\label{sec:SymplectoStuff}

The chamber complex for the multiplicity function can be used to identify domains of polynomiality.  However, these domains are not guaranteed to be as large as possible as seen in the examples of Section~\ref{sec:A2A3}.  In this section we improve the partition of the permutahedron into domains of polynomiality by identifying it as  a bounded plane arrangement that appears in symplectic geometry.  We begin by introducing the symplectic setup corresponding to the special case of type $A$ multiplicities.  Then we define the Duistermaat-Heckman function via an integral.  This function is piecewise polynomial with natural domains of polynomiality in terms of Weyl group orbits.  Finally, we will use a powerful theorem of Meinrenken \cite{Meinrenken} and Vergne \cite{Vergne}, the so-called \emph{Quantization Commutes with Reduction Theorem}, to show that the multiplicity function can be written locally as a very similar integral with the same domains.

Let $G=\mathrm{SU}(k)$, $T$ the Cartan subgroup of $G$, $\mathfrak{g}$ and $\mathfrak{t}$ their Lie algebras, $\mathfrak{t}_+^*$ the fundamental Weyl chamber and $\Lambda_W\subset\mathfrak{t}^*$ the weight lattice of $G$. For $\lambda\in\mathfrak{t}_+^*\cap\Lambda_W$, we will denote by $\Delta_\lambda$ the convex hull of the Weyl group orbit of $\lambda$ in $\mathfrak{t}^*$ (i.e. the permutahedron associated to $\lambda$). Let $O_\lambda = G\cdot \mathrm{diag}(\lambda)$ be the coadjoint orbit for $\lambda$.  We can view $O_\lambda$ as the set of $k\times k$ Hermitian matrices with eigenvalues $\{\lambda_1,\ldots,\lambda_k\}$.  By a theorem of Schur and Horn \cite{Schur,Horn} (or Kostant's convexity theorem \cite{Kostant2}, which extends the result to all compact Lie groups), $\Delta_\lambda$ is the image of the coadjoint orbit $O_\lambda$ with respect to the projection map 
\begin{equation}
\label{eqn:projection} \pi\,:\ \mathfrak{g}^* \longrightarrow
\mathfrak{t}^*\,.
\end{equation} 
The coadjoint orbits $O_\lambda$ are the geometric counterpart to the irreducible representations of $G$ with highest weight $\lambda$.  Note, the multiplicities for irreducible representations for $SU(k)$ and $SL(k)$ are the same.

Consider $M=O_\lambda$ and let $\Phi: M\rightarrow\mathfrak{t}^*$ be the restriction of $\pi$ to $M$.  In this case, $\Phi$ is the moment map of the symplectic manifold $M$ under the $T$ action.   The set $\Delta_{\mathrm{reg}}\subset\Delta_\lambda$ of regular values of $\Phi$ decomposes into a disjoint union of its connected components:
\begin{equation}
\Delta_{\mathrm{reg}} = \bigcup\Delta_i
\end{equation} 
and each $\Delta_i$ is an open convex polytope by a generalization of Kostant's convexity theorem due to Atiyah \cite{Atiyah} and Guillemin-Sternberg \cite{GS2}.  In fact, the singular values of $\Phi$ have the following nice combinatorial description.  This theorem first appeared in Heckman's thesis \cite{Hec}

\begin{theorem}{\rm\textbf{(\cite[Theorem 5.2.1]{GLS}, \cite{Hec})}}
\label{thm:regular.pts}
The singular points of the moment map $\Phi: M\rightarrow\mathfrak{t}^*$ are the convex polytopes 
\begin{equation}\label{eqn:walls}
\mathrm{conv}(W \cdot \sigma(\lambda) )
\end{equation}
where $\sigma \in \mathfrak{S}_k$ and $W$ is any parabolic subgroup of $\mathfrak{S}_k$ generated by all reflections corresponding to roots orthogonal to a conjugate of a fundamental weight.
\end{theorem}

In other words, the $\Delta_{i}$'s are the regions in the arrangement given by slicing the permutahedron by bounded hyperplane regions parallel to one of its exterior facets which pass through orbit points. See for example Figure~\ref{fig:A2diagrams} in Section~\ref{sec:A2A3}.  Note, this is not a hyperplane arrangement inside a polytope since the convex hulls do not necessarily extend to the boundary of the permutahedron.

Duistermaat and Heckman have shown that much of the geometry of coadjoint orbits can be determined simply by studying the $\Delta_{i}$'s.  For $\mu\in\Delta_i$, the \emph{symplectic reduction} of $M$ at the regular value $\mu$ of $\Phi$ is defined by 
\begin{equation}
M_\mu = \Phi^{-1}(\mu)/T \, .
\end{equation}
For arbitrary $G$, the reduced space $M_\mu$ of $M$ at a regular value of $\Phi$ is an orbifold, but for $\mathrm{SU}(k)$ this orbifold is a compact K\"ahler manifold whose symplectic form we will denote by $\omega_\mu$ \cite{DH1}.  Duistermaat and Heckman \cite{DH1} have shown that $M_\mu \cong M_{\mu_0}$ as complex manifolds for any pair $\mu_{0},\mu\in\Delta_i$.  Furthermore, they have also shown the linear variation formula \cite{DH1},
\begin{equation}
\label{eqn:linearvariation}
\omega_\mu = \omega_{\mu_0} + \langle\mu-\mu_0, c\rangle\,,
\end{equation}
where $c\in\mathfrak{t}\otimes\Omega^2(M_\mu)$ is the Chern form of the principal $T$-bundle $\Phi^{-1}(\mu)\rightarrow M_\mu$. Therefore, they use this fact about the symplectic forms to show that, for $M_\mu$ of dimension $2d$, the symplectic volume function
\begin{equation}\label{eqn:DHintegral}
f^{\mathrm{DH}}_\lambda(\mu) = \int_{M_{\mu_0}}\mathrm{exp}\,\omega_\mu = \int_{M_{\mu_0}}\frac{\omega_\mu^{d}}{d!}
\end{equation}
is a polynomial function on $\Delta_i$, called the \emph{Duistermaan-Heckman polynomial}. Note that the only aspect of this integral that depends specifically on $\mu$, and not just on which connected component of regular values contains it, is the symplectic form which is determined by \eqref{eqn:linearvariation}.  From the integral, one can show that the degrees of these polynomials are less than or equal to $(\mathrm{dim}\,M)/2-\mathrm{dim}\, G$.

Using a theory of quantization initiated by Kostant, Kirillov and Souriau (see \cite{Kostant3}, for instance), we can apply the same reasoning used by Duistermaat and Heckman to the multiplicity function.

\begin{theorem}
\label{thm:SamePartition} The partitions of the permutahedron for $\mathfrak{su}(k)$ (or $\sln{k}$) into its domains of polynomiality for the weight multiplicities and for the Duistermaat-Heckman measure are the same.  Namely, the domains are the connected components of regular points determined by \eqref{eqn:walls}.
\end{theorem}

\begin{proof}
Let $\mathrm{Td}(M_{\mu_0})$ be the Todd form of $M_{\mu_0}$. The Quantization Commutes with Reduction Theorem \cite{Meinrenken,Vergne} asserts that for $\mu\in\Delta_i$,
\begin{equation}
m_\lambda(\mu) = \int_{M_{\mu_0}}(\mathrm{exp}\,\omega_\mu)\,\mathrm{Td}(M_{\mu_0})\,.
\end{equation}
The right hand side is the Hirzebruch-Riemann-Roch number of $M_\mu$. The only factor in the integral which depends on $\mu$ is the symplectic form, everything else depends only on the region containing $\mu$.  Thus by \eqref{eqn:linearvariation} $m_\lambda(\mu)$ is a polynomial function of $\mu$ on $\Delta_i$ as with the Duistermaat-Heckman measure.
\end{proof}

\begin{remark}
\label{rem:optimal.partition}
This proof implies that the optimal domains of polynomiality for the multiplicity function must be unions of the $\Delta_{i}$'s.  Guillemin, Lerman and Sternberg \cite{GLS} have shown that this partition is optimal for the Duistermaat-Heckman measure by showing that the difference between the polynomials in two adjacent regions is nonzero.  We conjecture that this partition is also optimal for the multiplicity function.  This has been confirmed up to $\mathrm{SL}_4\mathbb{C}$.  
\end{remark}

As further evidence for the conjecture, we note that on a given domain, the weight polynomial and the Duistermaat-Heckman polynomial in \eqref{eqn:DHintegral} have the same leading term since 
\begin{equation}
\mathrm{Td}(M_\mu) = 1 + \sum_{j=1}^d\tau_j  
\end{equation}
with $\tau_j\in\Omega^{2j}(M_{\mu_0})$ in the de Rham complex.

\begin{remark}
\label{rem:counting.problem}
It is a very interesting open problem to count the regions in the permutahedron subdivided according to Theorem~\ref{thm:regular.pts}. This is the analog of Kirillov's question for the Kostant partition function mentioned in Section~\ref{subsec:KPF}.  We have determined all the region counts for $\mathrm{SL}_4\mathbb{C}$ in Figure~\ref{fig:number_regions_permutahedra}.
\end{remark}

There are many links between the weight multiplicities and the Duistermaat-Heckman function. For example, Dooley, Repka and Wildberger \cite{DRW} provide a way to go from the weight diagram for $\lambda$ to the Duistermaat-Heckman measure for $O_{\lambda+\delta}$:
\begin{equation}
f^{\mathrm{DH}}_{\lambda+\delta} = \sum_{\textrm{$\beta$ weight of $V_\lambda$}}m_\lambda(\beta)\,f^{\mathrm{DH}}_\delta\,.
\end{equation}
Also, if $\nu$ is the Lebesgue measure on $\mathfrak{t}^*$, normalized so that the parallelepiped given by the simple root vectors has unit measure, we define the \emph{Duistermaat-Heckman measure} to be the product $f^{\mathrm{DH}}\nu$. Now for each $n\in\mathbb{N}$ construct the discrete measure
\begin{equation}
\nu_n = \frac{1}{\mathrm{dim}\,V_{n\lambda}}\sum_{\textrm{$\beta$ weight of $V_{n\lambda}$}}m_{n\lambda}(\beta)\,\delta_{\beta/n}
\end{equation}
where $\delta_x$ is a point mass at $x$ and $V_{n\lambda}$ is the irreducible representation of $\mathfrak{su}(k)$ with highest weight $n\lambda$. Then Heckman \cite{Hec} proved that $\nu_n$ converges weakly to the Duistermaat-Heckman measure as $n\rightarrow\infty$.
 
Furthermore, the Duistermaat-Heckman function from above can be computed in the following way (see \cite{GLS}): its value at a point $(\lambda,\beta)$ is obtained by using Kostant's multiplicity formula with $\delta=0$ and a deformation of Kostant's partition function that takes the volume of the polytopes $\{(k_\alpha)_{\alpha\in\Delta_+}\in\mathbb{R}^{|\Delta_+|}_{\geq0}\ : \ \sum_{\alpha\in\Delta_+}k_\alpha\alpha = v\}$ instead of their number of integral points.

\section{The Kostant arrangements}
\label{sec:KostantArrangements}

In this section, we will construct a hyperplane arrangement whose regions are also domains of polynomiality for the multiplicity function. This partition into domains will be unlike the ones obtained in Remark~\ref{rem:ComplexPartition} and Theorem~\ref{thm:regular.pts} in that it is not invariant under rescaling $\lambda$ and $\beta$. We will deduce the form of this arrangement from a closer look at Kostant's multiplicity formula~(\ref{eqn:KMF}) and its chamber complex defined in Section~\ref{subsec:KPF}.  

\begin{lemma}
\label{lemma:KPFwalls}
The set of normals to the facets of the maximal cones of the chamber complex of the Kostant partition function of $A_n$ consists of all the conjugates of the fundamental weights.
\end{lemma}

\begin{proof}
The facets of the maximal cones of the chamber complex span the same hyperplanes as the facets of the base cones whose common refinement is the chamber complex. Base cones correspond to sets of $n$ linearly independent positive roots.  Fixing a particular base cone spanned by $\{\gamma_1,\ldots,\gamma_n\}$, consider the undirected graph $G$ on $\{1,\ldots,n+1\}$ where $(i,j)$ is an edge if $e_i-e_j=\gamma_m$ for some $m$. The fact that the $\gamma_j$'s are linearly independent implies that $G$ has no cycles. So $G$ is a forest, and since it has $n+1$ vertices and $n$ edges (one for each $\gamma_j$), it is actually a tree. Suppose now we remove $\gamma_j = e_s-e_t$ and want to find the normal of the hyperplane spanned by the other $\gamma_i$'s. The graph $G$ with the edge $(s,t)$ removed consists of two trees $T_1$ and $T_2$. List $\{1,\ldots,n+1\}$ in the form
\begin{displaymath}
\sigma : \ \underbrace{i_1, i_2, \ldots, i_{j-1}, s}_{\textrm{vertices of $T_1$}}, \underbrace{t, i_j, i_{j+1}, \ldots, i_{n+1-2}}_{\textrm{vertices of $T_2$}}
\end{displaymath}
where we will think of $\sigma$ as a permutation in one-line form. 

Now let $\alpha_i'=e_{\sigma(i)}-e_{\sigma(i+1)}$ and note that $\alpha_j'=e_s-e_t=\gamma_j$. The set $\{\alpha_1',\ldots,\alpha_n'\}$ is a root system basis because it is the image under the action of $\sigma^{-1}$ of the original simple roots $\alpha_i=e_i-e_{i+1}$. Observe that every edge in $T_1$ can be expressed as a sum of $\alpha_1',\ldots,\alpha_{j-1}'$, and every edge in $T_2$ as a sum of $\alpha_{j+1}',\ldots,\alpha_n'$, so that all $\gamma_i$'s in $\{\gamma_1,\ldots,\widehat{\gamma_j},\ldots,\gamma_n\}$ can be expressed as linear combinations of $\alpha_1',\ldots,\widehat{\alpha_j'},\ldots,\alpha_n'$. The normal for the corresponding hyperplane will therefore be the $j$th fundamental weight $\omega_j'$ for the basis $\{\alpha_1',\ldots,\alpha_n'\} = \sigma\cdot\{\alpha_1\ldots,\alpha_n\}$.

Conversely, given any fundamental weight $\omega_j'$ for the root system basis $\sigma\cdot\{\alpha_1\ldots,\alpha_n\}$ (or equivalently, $\sigma^{-1}\cdot\omega_j$, where $\omega_j$ is the $j$th fundamental weight for the standard simple roots), we want to show it can occur as the normal to a hyperplane. Let $H$ be a hyperplane separating the standard positive roots from the negative ones. For each $\alpha_i'=\sigma\cdot\alpha_i$, we can pick a sign $\varepsilon_i$ such that $\varepsilon_i\alpha_i'$ is on the positive side of $H$. Hence $\{\varepsilon_1\alpha_1',\ldots,\varepsilon_n\alpha_n'\}$ is a linearly independent subset of the set of standard positive roots, and thus it corresponds to one of the base cones of $M_{A_n}$. The corresponding graph is a path since we have a system of simple roots (up to sign reversal). Removing $\varepsilon_j\alpha_j'$ and applying the above procedure with the order given by the path gives that $\omega_j'$ occurs as the normal of the corresponding hyperplane.
\end{proof}

To compute multiplicities for $\sln{k}$ using Kostant's formula, we look at the points $\sigma(\lambda+\delta)-(\beta+\delta)$, as $\sigma$ ranges over the Weyl group $\mathfrak{S}_k$. Some of these points will lie inside the chamber complex for the Kostant partition function and we compute the multiplicity by finding which cells contain them and evaluating the corresponding polynomials at those points. Starting with generic $\lambda$ and $\beta$, none of the points $\sigma(\lambda+\delta)-(\beta+\delta)$ will lie on a wall of the chamber complex of the Kostant partition function, and if we move $\lambda$ and $\beta$ around a little in such a way that none of the $\sigma(\lambda+\delta)-(\beta+\delta)$ crosses a wall, we will obtain the multiplicity for the new $\lambda$ and $\beta$ by evaluating the same polynomials. So there is a neighborhood of $(\lambda,\beta)$ on which the multiplicity function is given by the same polynomial in variables $\lambda$ and $\beta$. 

Lemma~\ref{lemma:KPFwalls} describes the walls of the chamber complex for the Kostant partition function in terms of the normals to the hyperplanes (though the origin) supporting the facets of the maximal cells. Now a point $\sigma(\lambda+\delta)-(\beta+\delta)$ will be on one of those walls (hyperplane though the origin) when its scalar product with the hyperplane's normal, say $\theta(\omega_j)$, vanishes, that is when
\begin{equation}\label{eqn:KosArrWalls}
\langle \sigma(\lambda+\delta)-(\beta+\delta), \theta(\omega_j)\rangle = 0
\end{equation}

\begin{figure}[!ht]
\begin{center}
\begin{tabular}{c@{\hspace{10mm}}c}
\fboxsep10pt\fbox{\includegraphics[height=0.35\textwidth]{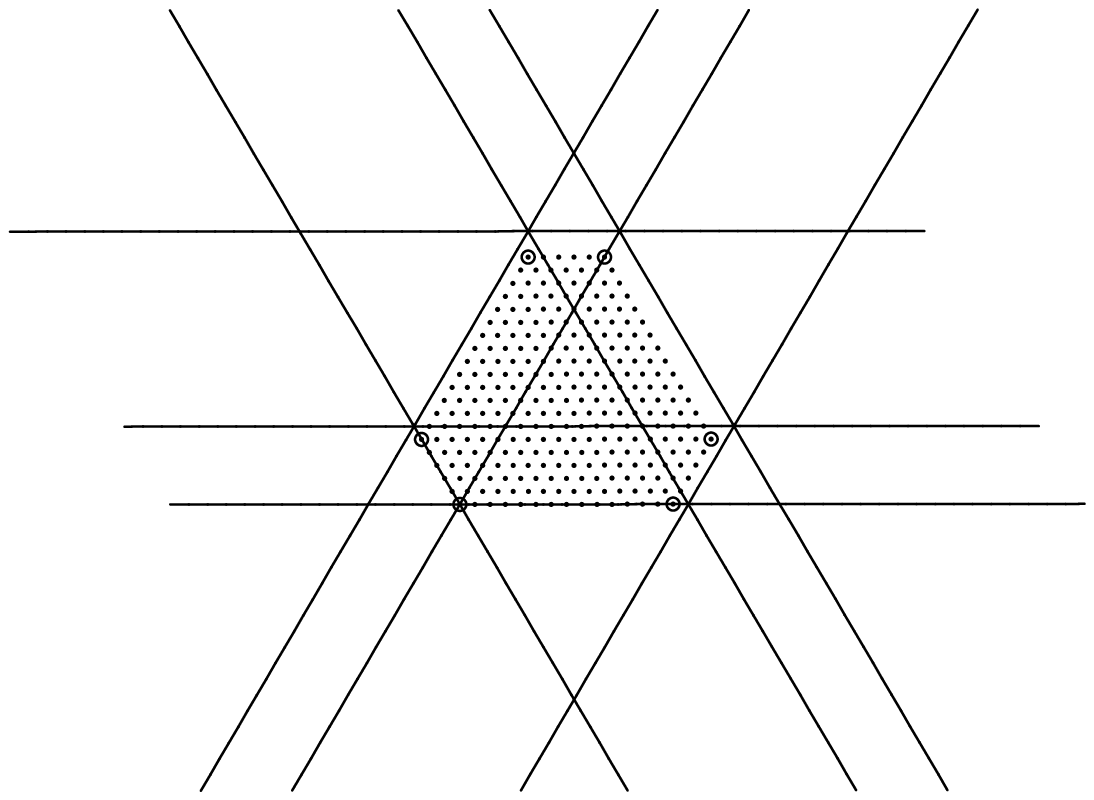}} & \fboxsep10pt\fbox{\includegraphics[height=0.35\textwidth]{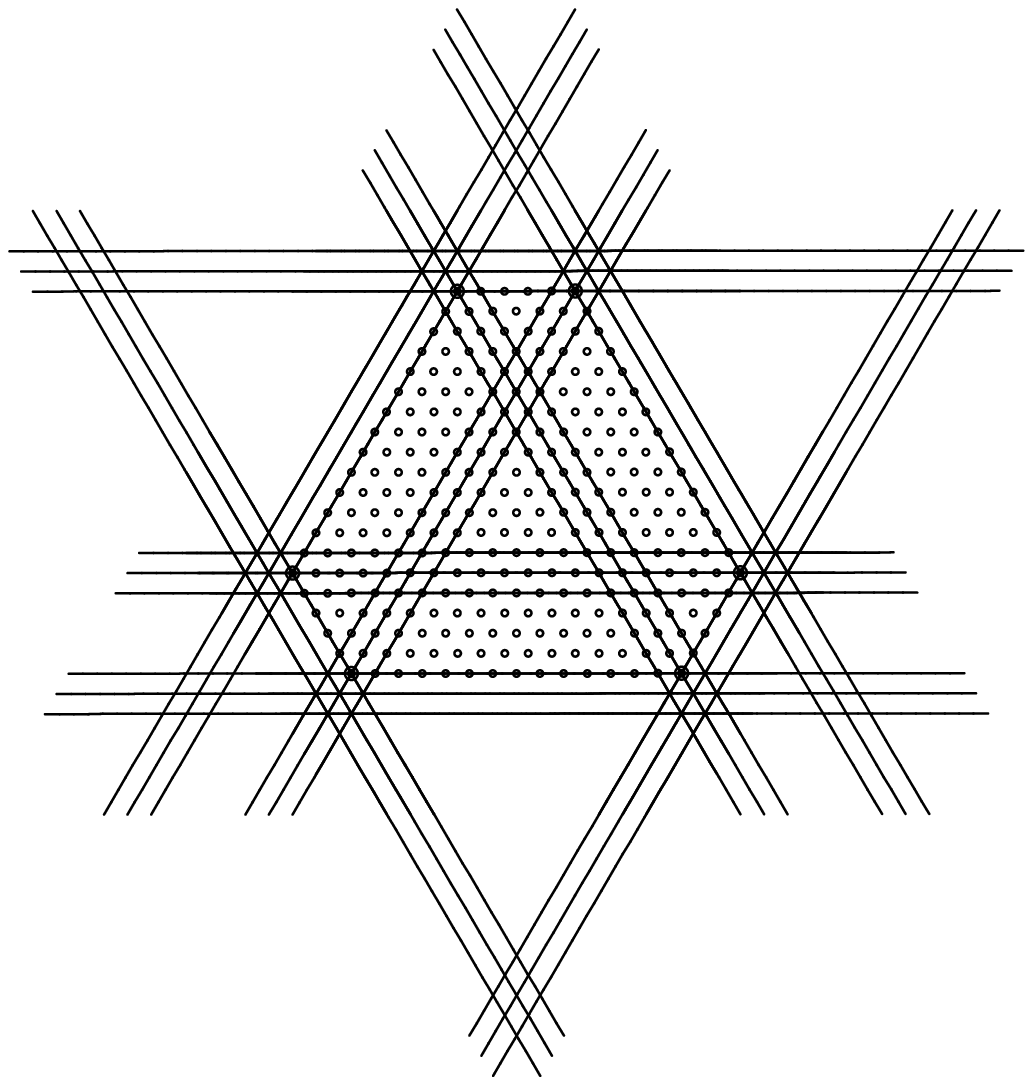}} \\
\end{tabular}
\caption{Kostant arrangement $\mathcal{A}^{(id)}_{(11,-3,-8)}$ for $A_2$ (left). Superposition of the Kostant arrangements $\mathcal{A}^{(\psi)}_{(11,-3,-8)}$ for all choices of $\psi$ (right).}
\label{fig:KosArr1}
\end{center}
\end{figure}

For any $\lambda$, consider the arrangement of all such hyperplanes for $1\leq j\leq k$ and $\sigma,\,\theta\in\mathfrak{S}_k$. For $\beta$ and $\beta'$ in the same region of this arrangement and any fixed $\sigma\in\mathfrak{S}_k$, the points $\sigma(\lambda+\delta)-(\beta+\delta)$ and $\sigma(\lambda+\delta)-(\beta'+\delta)$ lie on the same side of every wall of the chamber complex for the Kostant partition function. Figure~\ref{fig:KosArr1} (on the left) shows the arrangement we get for $\lambda=(11,-3,-8)$ in $A_2$, with and without the weight diagram.

In view of the invariance of the multiplicities under the action of the Weyl group, Kostant's formula has to give the same thing if we replace $\beta$ by $\psi(\beta)$, $\psi\in\mathfrak{S}_k$. Replacing $\beta$ by $\psi(\beta)$ in the equations~(\ref{eqn:KosArrWalls}) above yields another hyperplane arrangement, which we will denote by $\mathcal{A}^{(\psi)}_\lambda$. Hence for each $\lambda$, we get a family of arrangements indexed by $\psi\in\mathfrak{S}_k$, which we will call the \emph{Kostant arrangements} for $\lambda$. Figure~\ref{fig:KosArr1} (on the right) shows the superposition of those arrangements for $\lambda=(11,-3,-8)$ as $\psi$ ranges over the Weyl group.

Suppose the chamber complex for the Kostant partition function for $\sln{k}$ has $r(k)$ full dimensional cones. We will choose a labeling of these regions with the integers $1, \ldots, r(k)$ once and for all, and let the associated polynomials be $p_1,\ldots,p_{r(k)}$. Recall, these are polynomials of degree ${k-1 \choose 2}$ on the subspace $x_1+\cdots+x_k=0$ of $\mathbb{R}^k$. We will also label the exterior of the chamber complex by $0$ and let its polynomial be the zero polynomial $p_0$.

\begin{definition}
For generic $\lambda$ and $\beta$, let $v^{(\psi)}_\sigma(\lambda,\beta)$ (or just $v^{(\psi)}_\sigma$) be the label of the region containing the point $\sigma(\lambda+\delta)-(\psi(\beta)+\delta)$ (this label is unique for generic $\lambda$ and $\beta$). Define the \emph{type} of $\lambda$ and $\beta$ to be the vector
\begin{displaymath}
\mathrm{Type}^{(\psi)}(\lambda,\beta) = (v^{(\psi)}_\sigma)_{\sigma\in\mathfrak{S}_k}\,,
\end{displaymath}
for some fixed total order on $\mathfrak{S}_k$. Furthermore, define 
\begin{equation}
\label{eqn:TypePolynomial}
P^{(\psi)}_\lambda(\beta) = \sum_{\sigma\in\mathfrak{S}_k}(-1)^{l(\sigma)}p_{v^{(\psi)}_\sigma}(\sigma(\lambda+\delta)-(\psi(\beta)+\delta))\,.
\end{equation}
\end{definition}

\begin{lemma}
\label{lemma:KostantPolynomial}
$P^{(\psi)}_\lambda$ is a polynomial function on the interior of the regions of $\mathcal{A}^{(\psi)}_\lambda$ and coincides with the multiplicity function there.
\end{lemma}

\begin{proof}
For fixed $\lambda$, the type of points along a path between two $\beta$'s in the interior of the same region of $\mathcal{A}^{(\psi)}_\lambda$ will remain the same by definition of the Kostant arrangement (because no $\sigma(\lambda+\delta)-(\psi(\beta)+\delta)$ crosses a wall along that path).
\end{proof}

The reason why Lemma~\ref{lemma:KostantPolynomial} is restricted to the interior of the regions is that while polynomials for adjacent regions of the chamber complex for the Kostant partition function have to coincide on the intersection of their closures, there is a discontinuous jump in the value of the Kostant partition function (as a piecewise polynomial function) when going from a region on the boundary of the complex to region $0$ (outside the complex). 

\begin{remark}
\label{rem:EhrhartGT}
Given a rational polytope $Q$ of dimension $d$ in $\mathbb{R}^n$ and $t\in\mathbb{N}$, denote by $tQ$ the polytope obtained by scaling $Q$ by a factor of $t$. Ehrhart \cite{Ehrhart} showed that the function $t\mapsto|tQ\cap\mathbb{Z}^n|$, counting the number of integer points in $tQ$ as a function of $t$, is a quasipolynomial of degree $d$, and a polynomial of degree $d$ if $Q$ is integral. This function is called the \emph{Ehrhart (quasi)polynomial} of the polytope $Q$. Furthermore, the leading coefficient of the Ehrhart quasipolynomial is the $d$-dimensional volume of $Q$. It can be shown that for every fixed $\lambda$ and $\beta$, and any $\psi$, the function $t\mapsto\mathrm{Type}^{(\psi)}(t\lambda,t\beta)$ in the nonnegative integer variable $t$ eventually stabilizes as $t$ grows (in a way that depends only on $k$ and not on $\lambda$ and $\beta$). This can be used to give a proof that the Ehrhart functions of the Gelfand-Tsetlin polytopes $\mathrm{GT}_{\lambda,\beta}$ are polynomial (we omit the proof since we prove something stronger in Corollary~\ref{cor:scaling} below).
\end{remark}

In the definition of the Kostant arrangements above, a lot of the hyperplanes are redundant. We simplify here the description of these arrangements. Since everything occurs on the subspace $x_1+\cdots+x_k=0$ of $\mathbb{R}^k$, so that $\langle \sigma(\lambda+\delta)-(\psi(\beta)+\delta), (1,1,\ldots,1)\rangle = 0$, we can regard the normals to the hyperplanes up to adding multiples of $(1,1,\ldots,1)$ without changing the Kostant arrangements. So we can use $\tilde{\omega}_j = \omega_j+\frac{j}{k}(1,1,\ldots,1)$:
\begin{displaymath}
\tilde{\omega}_j = \frac{1}{k}(\underbrace{k,k,\ldots,k}_{\textrm{$j$ times}},\underbrace{0,0,\ldots,0}_{\textrm{$k-j$ times}}) = (\underbrace{1,1,\ldots,1}_{\textrm{$j$ times}},\underbrace{0,0,\ldots,0}_{\textrm{$k-j$ times}}) = e_1+\cdots+e_j\,,
\end{displaymath}
which is more convenient than $\omega_j$ for what follows. The hyperplanes of $\mathcal{A}^{(\psi)}_\lambda$ then have the form
\begin{eqnarray}
\label{eqn:KosHypPsi}
0 & = & \langle \sigma(\lambda+\delta)-(\psi(\beta)+\delta), \theta(\tilde{\omega}_j)\rangle\nonumber\\ 
0 & = & \langle \sigma(\lambda)-\psi(\beta)+\sigma(\delta)-\delta, e_{\theta(1)}+\cdots+e_{\theta(j)}\rangle\nonumber\\
0 & = & \langle (\lambda_{\sigma^{-1}(i)}-\beta_{\psi^{-1}(i)}+\delta_{\sigma^{-1}(i)}-\delta_i)_{i=1,\ldots,k}, e_{\theta(1)}+\cdots+e_{\theta(j)}\rangle\nonumber\\ 
0 & = & \sum_{i=1}^j(\lambda_{\sigma^{-1}(\theta(i))}-\beta_{\psi^{-1}(\theta(i))}+\delta_{\sigma^{-1}(\theta(i))}-\delta_{\theta(i)})\nonumber\\
\beta_{\psi^{-1}(\theta(1))}+\cdots+\beta_{\psi^{-1}(\theta(j))} & = & \lambda_{\sigma^{-1}(\theta(1))}+\cdots+\lambda_{\sigma^{-1}(\theta(j))} + \sum_{i=1}^j(\delta_{\sigma^{-1}(\theta(i))}-\delta_{\theta(i)})\,.
\end{eqnarray}

At this point we can get rid of the permutations since only the subsets $\theta(\{1,\ldots,j\})$, $\psi^{-1}\theta(\{1,\ldots,j\})$ and $\sigma^{-1}\theta(\{1,\ldots,j\})$ are important and not the order of their elements. They can be any subsets since $\psi^{-1}\theta$, $\sigma^{-1}\theta$ and $\theta$ can be any three permutations of $\mathfrak{S}_k$. Because the $\beta_i$, the $\lambda_i$ and the $\delta_i$ sum up to zero, replacing these subsets by their complements gives the same hyperplane.  This proves the following proposition.

\begin{proposition}
\label{prop:KosHyp}
The hyperplanes of the Kostant arrangements are defined by the equations 
\begin{equation}
\label{eqn:KosHyp}
\beta_{u_1}+\cdots+\beta_{u_j} = \lambda_{v_1}+\cdots+\lambda_{v_j} + \sum_{i=1}^j(\delta_{v_i}-\delta_{w_i})\,,
\end{equation}
where $U=\{u_1,\ldots,u_j\}$, $V=\{v_1,\ldots,v_j\}$ and $W=\{w_1,\ldots,w_j\}$ range over all $j$-element subsets of $\{1,\ldots,k\}$ and  $j \leq \lfloor k/2\rfloor$.
\end{proposition}

We will call the correction term involving only $\delta$ in a hyperplane given as in~(\ref{eqn:KosHyp}), the $\emph{$\delta$-shift}$:
\begin{equation}
\mathrm{shift}(V,W) = \sum_{i=1}^j(\delta_{v_i}-\delta_{w_i})\,.
\end{equation}

\begin{remark}
\label{rem:outershift}
For fixed $U$, we get a series of parallel hyperplanes, and we can determine which are the outer ones because they correspond to maximal and minimal sums of $\lambda_{q_i}$. Since the $\lambda_1\geq\lambda_2\geq\cdots\geq\lambda_k$, they are 
\begin{equation}
\label{eqn:OuterHyp}
\begin{array}{rcl}
\beta_{u_1}+\cdots+\beta_{u_j} = & \lambda_1+\cdots+\lambda_j & + \quad \mathrm{shift}(\{1,\ldots,j\},W)\\[2mm]
\beta_{u_1}+\cdots+\beta_{u_j} = & \lambda_{k-j+1}+\cdots+\lambda_k & + \quad \mathrm{shift}(\{k-j+1,\ldots,k\},W)\,.
\end{array}
\end{equation}

Note that since the coordinates of $\delta$ are decreasing, $\mathrm{shift}(\{1,\ldots,j\},W)\geq 0$ \ and \ $\mathrm{shift}(\{k-j+1,\ldots,k\},W)\leq 0$ for all $W$. 
\end{remark} 

We conclude this section by relating the domains given by the Kostant arrangements and those given by Theorem~\ref{thm:SamePartition}, by showing that the hyperplanes supporting the facets of the domains are precisely the hyperplanes of the Kostant arrangements without the $\delta$-shift factors.

\begin{proposition}
\label{prop:SupportingHyperplanes}
The supporting hyperplanes of the facets of the top-dimensional domains of the permutahedron for generic $\lambda$ are the hyperplanes
\begin{equation}
\beta_{u_1}+\cdots+\beta_{u_j} = \lambda_{v_1}+\cdots+\lambda_{v_j}\,, 
\end{equation}
for $1\leq j\leq\lfloor k/2\rfloor$ and $U=\{u_1,\ldots,u_j\}$, $V=\{v_1,\ldots,v_j\}$ ranging over all pairs of $j$-element subsets of $\{1,\ldots,k\}$.
\end{proposition}

\begin{proof}
Theorem~\ref{thm:regular.pts} gives the walls supporting the facets as the convex hulls of $W\cdot\sigma(\lambda)$, where $\sigma(\lambda)$ is a point of the Weyl orbit of $\lambda$, and $W$ is a parabolic subgroup of the Weyl group. For $\mathfrak{S}_k$, those subgroups permute two complementary sets of indices independently. If $U$ is one of those sets of indices, with $|U|=j$, and $\lambda_{v_1},\ldots,\lambda_{v_j}$ the coordinates of $\sigma(\lambda)$ in those positions, then the hyperplane supporting $W\cdot\sigma(\lambda)$ is
\begin{equation}
\beta_{u_1}+\cdots+\beta_{u_j} = \lambda_{v_1}+\cdots+\lambda_{v_j}\,. 
\end{equation}
Had we chosen the complement of $U$ instead with the remaining $\lambda_i$'s, we would have gotten the same hyperplane in the subspace $x_1+\cdots+x_k=0$ of $\mathbb{R}^k$ since the $\lambda_i$'s and the $\beta_i$'s sum up to zero.
\end{proof}

We can obtain the following corollary without using the full description of the domains of the permutahedron obtained by symplectic geometry means in Theorem~\ref{thm:regular.pts}. 

\begin{corollary}
\label{cor:BoundaryWalls}
The hyperplanes supporting the facets of the permutahedron for a generic $\lambda$ are 
\begin{equation}
\begin{array}{rcl}
\beta_{u_1}+\cdots+\beta_{u_j} & = & \lambda_1+\cdots+\lambda_j \\
\beta_{u_1}+\cdots+\beta_{u_j} & = & \lambda_{k-j+1}+\cdots+\lambda_k
\end{array}
\end{equation}
for $1\leq j\leq\lfloor k/2\rfloor$ and $U=\{u_1,\ldots,u_j\}$ ranging over all $j$-element subsets of $\{1,\ldots,k\}$.
\end{corollary}

\begin{proof}
We remark that the ``shell'' of the weight diagram is just a permutahedron, whose facets can easily be described in terms of permutations (see \cite[p.~18]{Ziegler}). For $U\subseteq\{1,\ldots,k\}$, construct a $k$-vector by putting the first $|U|$ $\lambda_i$'s in the positions indexed by $U$ and filling the other positions with the remaining elements. Then act by the subgroup of $\mathfrak{S}_k$ that permutes the elements in positions $U$ and $\{1,\ldots,k\}\setminus U$ independently to get a facet as the convex hull of the points of this orbit. The affine span of this facet is the hyperplane
\begin{displaymath}
\beta_{u_1}+\cdots+\beta_{u_j} = \lambda_1+\cdots+\lambda_j\,.
\end{displaymath}

By choosing the last $|U|$ $\lambda_i$'s instead, we get the hyperplane supporting the opposite parallel facet. These are the outer hyperplanes from~(\ref{eqn:OuterHyp}) without the shifts. Remark~\ref{rem:outershift} also implies these outer hyperplanes actually lie outside the permutahedron.
\end{proof}

\begin{remark}
\label{rem:directions}
The Weyl orbits of $\omega_j$ and $\omega_{k-j}$ ($1\leq j\leq k-1$) for $\sln{k}$ (type $A_{k-1}$) determine the same set of directions, since $\omega_{k-j}$ is $-\omega_j$ with the coordinates in reverse order. So the Weyl orbits of $\omega_1,\ldots,\omega_{\lfloor k/2\rfloor}$ already determine all the possible normals to facets of the permutahedron (and the hyperplanes of the Kostant arrangement).
\end{remark}

\section{Polynomiality in the chamber complex}
\label{sec:PolynomialityComplex}

Theorem~\ref{thm:SPF} allowed us to write the multiplicity function as a partition function, which is therefore quasipolynomial over the convex polyhedral cones of the chamber complex $\mathcal{C}^{(k)}$. On the other hand, for each dominant weight $\lambda$, Theorem~\ref{thm:SamePartition} shows that the partition of the permutahedron from Theorem~\ref{thm:regular.pts} gives domains over which the multiplicity function is polynomial in $\beta$. We show here that the quasipolynomials attached to the complex $\mathcal{C}^{(k)}$ are actually polynomials, so that the multiplicity function is polynomial in both $\lambda$ and $\beta$ over the cones of the complex. 

The union of the cones of the complex $\mathcal{C}^{(k)}$ is the cone
\begin{equation}
T^{(k)} = \bigcup_{\lambda\in C_0}\{\lambda\}\times\mathrm{conv}(\mathfrak{S}_k\cdot\lambda)\,,
\end{equation}
where $C_0$ is the fundamental Weyl chamber.

We can lift the partition of the permutahedron from Theorem~\ref{thm:regular.pts} to $(\lambda,\beta)$-space by lifting the wall
\begin{displaymath}
\mathrm{conv}(W\cdot\sigma(\lambda))
\end{displaymath} 
to
\begin{displaymath}
\bigcup_{\lambda\in C_0}\{\lambda\}\times\mathrm{conv}(W\cdot\sigma(\lambda))\,.
\end{displaymath}
This gives a partition $\mathcal{T}^{(k)}$ of the cone $T^{(k)}$ into convex polyhedral cones, and Theorem~\ref{thm:SPF} implies that the multiplicity function is quasipolynomial over the cones of $\mathcal{T}^{(k)}$. We recover the domains from Theorem~\ref{thm:regular.pts} by intersecting $\mathcal{T}^{(k)}$ with $L(\lambda)$ from equation~\eqref{eqn:Llambda}. Our reason for introducting $\mathcal{T}^{(k)}$ rather that working with the complex $\mathcal{C}^{(k)}$ is that Proposition~\ref{prop:SupportingHyperplanes} lets us describe the hyperplanes supporting the facets of the cones of $\mathcal{T}^{(k)}$ easily. Indeed, if 
\begin{displaymath}
\beta_{u_1}+\cdots+\beta_{u_j} = \lambda_{v_1}+\cdots+\lambda_{v_j}
\end{displaymath}
supports a wall $\mathrm{conv}(W\cdot\sigma(\lambda))$ for fixed $\lambda$, then the wall $\bigcup_{\lambda\in C_0}\{\lambda\}\times\mathrm{conv}(W\cdot\sigma(\lambda))$ is supported by the hyperplane
\begin{displaymath}
\beta_{u_1}+\cdots+\beta_{u_j} = \lambda_{v_1}+\cdots+\lambda_{v_j}
\end{displaymath}
in $(\lambda,\beta)$-space, where we now think of $\lambda$ as variable, just like $\beta$.

The last tool we need is a lifted version of the Kostant arrangements. Recall from~\eqref{eqn:KosHypPsi} that the Kostant arrangement $\mathcal{A}^{(\psi)}_\lambda$ has the hyperplanes
\begin{displaymath}
\beta_{\psi^{-1}(\theta(1))}+\cdots+\beta_{\psi^{-1}(\theta(j))} = \lambda_{\sigma^{-1}(\theta(1))}+\cdots+\lambda_{\sigma^{-1}(\theta(j))} + \sum_{i=1}^j(\delta_{\sigma^{-1}(\theta(i))}-\delta_{\theta(i)})\,,
\end{displaymath}
as $\theta$ ranges over $\mathfrak{S}_k$. We will denote by $\mathcal{A}^{(\psi)}$ the arrangement with hyperplanes
\begin{displaymath}
\beta_{\psi^{-1}(\theta(1))}+\cdots+\beta_{\psi^{-1}(\theta(j))} = \lambda_{\sigma^{-1}(\theta(1))}+\cdots+\lambda_{\sigma^{-1}(\theta(j))} + \sum_{i=1}^j(\delta_{\sigma^{-1}(\theta(i))}-\delta_{\theta(i)})\,,
\end{displaymath}
where we now think of $\lambda$ as variable, and $\theta$ ranges over $\mathfrak{S}_k$ as before. The definition of $P^{(\psi)}_\lambda$ (equation~\eqref{eqn:TypePolynomial}) and Lemma~\ref{lemma:KostantPolynomial} generalize to give us a piecewise polynomial function $P^{(\psi)}$ in $\lambda$ and $\beta$ that expresses the multiplicity function as a polynomial on the interior of the regions of the arrangement $\mathcal{A}^{(\psi)}$.

\begin{theorem}
\label{thm:PolynomialityComplex}
The quasipolynomials determining the multiplicity function in the cones of $\mathcal{T}^{(k)}$ and $\mathcal{C}^{(k)}$ are polynomials of degree ${k-1 \choose 2}$ in the $\beta_i$, with coefficients of degree ${k-1 \choose 2}$ in the $\lambda_j$.
\end{theorem}

\begin{proof}
We will show that for each cone $C$ of $\mathcal{T}^{(k)}$ we can find a region $R$ of the Kostant arrangement $\mathcal{A}^{(\psi)}$ (for any $\psi$), such that $C\cap R$ contains an arbitrarily large ball. Then $P^{(\psi)}$ and the quasipolynomial in $C$ agree on the points $(\lambda,\beta)$ in that ball for which $(\lambda,\beta)\in\Lambda_W\times\Lambda_W$ and $\lambda-\beta\in\Lambda_R$ (the points corresponding to allowable pairs of a dominant weight and a weight of its irreducible representation). The quasipolynomial must therefore be polynomial on those points. The degree bounds follow from Remark~\ref{rem:LocallyPolynomial}.

By the remarks preceding this theorem, the hyperplanes supporting the facets of the cones of $\mathcal{T}^{(k)}$ are exactly the same as the hyperplanes of the Kostant arrangement $\mathcal{A}^{(\psi)}$ with the shifts removed. If we deform $\mathcal{A}^{(\psi)}$ continuously to make the shifts zero (by multipliying them by $t$ and letting $t$ going from 1 to 0, for example), the final deformed arrangement is a partition of $T^{(k)}$ that refines $\mathcal{T}^{(k)}$. Let $R$ be any region of $\mathcal{A}^{(\psi)}$ whose deformed final version is contained in $C$. Consider a ball of radius $r$ inside the deformed image of $R$, and suppose it is centered at the point $x$. If $s$ is the maximal amount by which the hyperplanes of the Kostant arrangement are shifted, then $R$ contains the ball of radius $r-s$ centered at $x$, and so does $C\cap R$. Since $C$ is a cone, we can make $r$ arbitrary large and the result follows since $s$ is bounded.

We get the same result for the complex $\mathcal{C}^{(k)}$ by passing to its common refinement with $\mathcal{T}^{(k)}$.
\end{proof}

Recall from Section~\ref{sec:SymplectoStuff} that the weight multiplicity function and the Duistermaat-Heckman function have the same leading term. In particular, the degree of the multiplicity function is at most the upper bound on the degree of the Duistermaat-Heckman function. For a torus $T$ acting on a symplectic manifold $M$, the latter is known to be $(\mathrm{dim}\,M)/2 - \mathrm{dim}\,T$. In our case, $M$ is the coadjoint orbit $O_\lambda$ and $\mathrm{dim}\,T = k-1$ since $T$ is the set of $k\times k$ traceless diagonal Hermitian matrices. The dimension of $O_\lambda$ is $k^2-k = k(k-1)$ for generic $\lambda$, but for nongeneric $\lambda$, we can get more precise bounds on the degrees. Since the coordinates of $\lambda$ are decreasing, it has the form
\begin{displaymath}
(\underbrace{\nu_1,\ldots,\nu_1}_{\textrm{$k_1$ times}},\underbrace{\nu_2,\ldots,\nu_2}_{\textrm{$k_2$ times}},\ldots,\underbrace{\nu_l,\ldots,\nu_l}_{\textrm{$k_l$ times}})
\end{displaymath}
where $\nu_1>\nu_2>\cdots>\nu_l$ and the $k_j$ sum up to $k$. In this case, one can show that $\mathrm{dim}\,O_\lambda = k^2-\sum k_j^2$, so that the weight multiplicity function for that $\lambda$ is piecewise polynomial of degree at most
\begin{equation}
\frac{k^2-\sum k_j^2}{2} - k + 1\,.
\end{equation}
For $\sln{4}$, for example, we get at most cubic polynomials for generic $\lambda$, at most quadratic polynomials for $\lambda$ with exactly two equal coordinates, at most linear polynomials for $\lambda$ with two pairs of equal coordinates ($\lambda$ of the form $(\nu,\nu,-\nu,-\nu)$) and constant polynomials for $\lambda$ with three equal coordinates ($\lambda$ of the form $(3\nu,-\nu,-\nu,-\nu)$ or $(\nu,\nu,\nu,-3\nu)$).

We can also deduce from Theorems~\ref{thm:SPF} and~\ref{thm:PolynomialityComplex} that the multiplicity function for type $A$ exhibits a scaling property in the following sense.

\begin{corollary}
\label{cor:scaling}
Let $\Upsilon$ be the set $\{(\lambda,\beta)\in\Lambda_W^2\ : \ \lambda-\beta\in\Lambda_R\}$. For any generic $(\lambda,\beta)\in\Upsilon$, we can find a neighborhood $U$\! of that point over which the function 
\begin{equation}
(\lambda,\beta,t)\in (U\cap\Upsilon)\times\mathbb{N} \longmapsto m_{t\lambda}(t\beta)
\end{equation}
is polynomial of degree at most $2{k-1\choose 2}$ in $t$ and ${k-1\choose 2}$ in the $\lambda$ and $\beta$ coordinates.
\end{corollary}

\begin{proof}
Let $(\lambda,\beta) \in \Upsilon$. For $U$ sufficiently small, the points $\{(t\lambda,t\beta)\ : \ t\in\mathbb{N}\}$ lie in the same cone of the chamber complex $\mathcal{C}^{(k)}$, and for $t\in\mathbb{N}$, $t\lambda$ and $t\beta$ are points on the weight lattice with their difference on the root lattice. Hence the corresponding multiplicities are obtained by evaluating the same polynomial at those points.  
\end{proof}

\begin{remark}
This corollary implies in particular that the Ehrhart functions (see Remark~\ref{rem:EhrhartGT}) of the Gelfand-Tsetlin polytopes $\mathrm{GT}_{\lambda,\beta}$ are always polynomial, even though the polytopes are not always integral (see \cite{DLMA}). 
\end{remark}

\section{Factorizations of weight polynomials}
\label{sec:Factorization}

In this section we use the explicit relation between the hyperplanes of the Kostant arrangements and the supporting hyperplanes of the partitioned permutahedron to identify certain factors in the weight polynomials. As mentioned in the introduction, Szenes and Vergne \cite{SzenesVergne} have recently observed this factorization phenomenon for general partition functions. The quasipolynomials associated to the partition function's chamber complex exhibit a certain number of linear factors that vanish on hyperplanes parallel and close to those supporting the walls of the complex. In our case it is unclear how to deduce the form of the walls of the complex $\mathcal{C}^{(k)}$ from the complex of the partition function given by matrix $E_k$ in Section~\ref{sec:SPF}. We are however able to deduce similar results from the Kostant arrangements and the description of the hyperplanes supporting the walls partitioning the permutahedron from Section~\ref{sec:KostantArrangements}.

\subsection{On the boundary of the permutahedron}

We have seen in Proposition~\ref{cor:BoundaryWalls} that each facet of the permutahedron is parallel and close to a hyperplane of a Kostant arrangement. This means that the domains of polynomiality of the weight diagram that are on the boundary of the permutahedron overlap with regions of the Kostant arrangement, but can't coincide because of the shifts caused by $\delta$. We can use this to our advantage to show that those weight polynomials have to factor somewhat. The reason is that two polynomials give the weights in the overlap: the one attached to a cone of the chamber complex obtained from writing the multiplicity function as a single partition function, and one, $P^{(\psi)}$, coming from Kostant's multiplicity formula. Because the overlap isn't perfect, the polynomial from Kostant's formula is valid on a region that goes outside the weight diagram and must therefore vanish there. The purpose of this section is to make precise this phenomenon and quantify it.

\begin{definition}
For fixed $\lambda$, consider the hyperplane
\begin{displaymath}
H\ : \ \beta_{u_1}+\cdots+\beta_{u_j} = \lambda_{v_1}+\cdots+\lambda_{v_j}
\end{displaymath}
where $1\leq j\leq\lfloor k/2\rfloor$ and $U=\{u_1,\ldots,u_j\}$, $V=\{v_1,\ldots,v_j\}$ are $j$-element subsets of $\{1,\ldots,k\}$. We will call the polynomial
\begin{equation}
\gamma_{U,V}(\lambda) = \beta_{u_1}+\cdots+\beta_{u_j} - \lambda_{v_1}-\cdots-\lambda_{v_j} \in \mathbb{Z}[\beta]
\end{equation} 
the \emph{defining equation} of $H$. For variable $\lambda$, we also define
\begin{equation}
\gamma_{U,V} = \beta_{u_1}+\cdots+\beta_{u_j} - \lambda_{v_1}-\cdots-\lambda_{v_j} \in \mathbb{Z}[\lambda,\beta]\,.
\end{equation}
\end{definition}

\begin{theorem}
\label{thm:factorization}
Let $R$ be a domain of polynomiality for the weight diagram of the irreducible representation of $\sln{k}$ with highest weight $\lambda$, and $p_R$ be its weight polynomial. Suppose that $R$ has a facet lying on the boundary of the permutahedron for $\lambda$ that has $\theta(\omega_j)$ as its normal vector, for some $\theta\in\mathfrak{S}_k$. If $\gamma = \gamma(\lambda)$ is the defining equation of the hyperplane supporting that facet,  then $p_R$ is divisible by the $j(k-j)-1$ linear factors $\gamma+1$, $\gamma+2$, $\ldots$, $\gamma+j(k-j)-1$, or $\gamma-1$, $\gamma-2$, $\ldots$, $\gamma-j(k-j)+1$.
\end{theorem}

Observe that this is invariant under replacing $j$ by $k-j$, which is a consequence of the remark in Remark~\ref{rem:directions}. By that remark, we can therefore restrict ourselves to $1\leq j\leq\lfloor k/2\rfloor$.

\begin{proof}
Suppose the hyperplane supporting the facet $F$ of $R$ on the boundary of the permutahedron has normal $\theta(\omega_j)$. By Remark~\ref{rem:outershift} and Proposition~\ref{cor:BoundaryWalls}, this hyperplane is either
\begin{eqnarray}\label{eqn:hyp1}
\beta_{\theta(1)}+\cdots+\beta_{\theta(j)} & = & \lambda_1+\cdots+\lambda_j \qquad\textrm{or}\\
\beta_{\theta(1)}+\cdots+\beta_{\theta(j)} & = & \lambda_{k-j+1}+\cdots+\lambda_k\,.\label{eqn:hyp2}
\end{eqnarray}

Suppose it's the first one (the argument is the same for the second one). From Proposition~\ref{prop:KosHyp}, we know that in the Kostant arrangements, we have the hyperplanes
\begin{equation}
\label{eqn:layers}
\beta_{\theta(1)}+\cdots+\beta_{\theta(j)} = \lambda_1+\cdots+\lambda_j + \mathrm{shift}(\{1,\ldots,j\}, W)\,,
\end{equation}
for $W$ ranging over $j$-element subsets of $\{1,\ldots,k\}$. We want to identify a region $R'$ of one of the Kostant arrangements that overlaps with $R$ and extends beyond the boundary of the weight diagram as far as possible. Note that although the exterior walls of $R$ and $R'$ have to be parallel, the interior walls don't. 

\begin{center}
\includegraphics[width=0.3\textwidth]{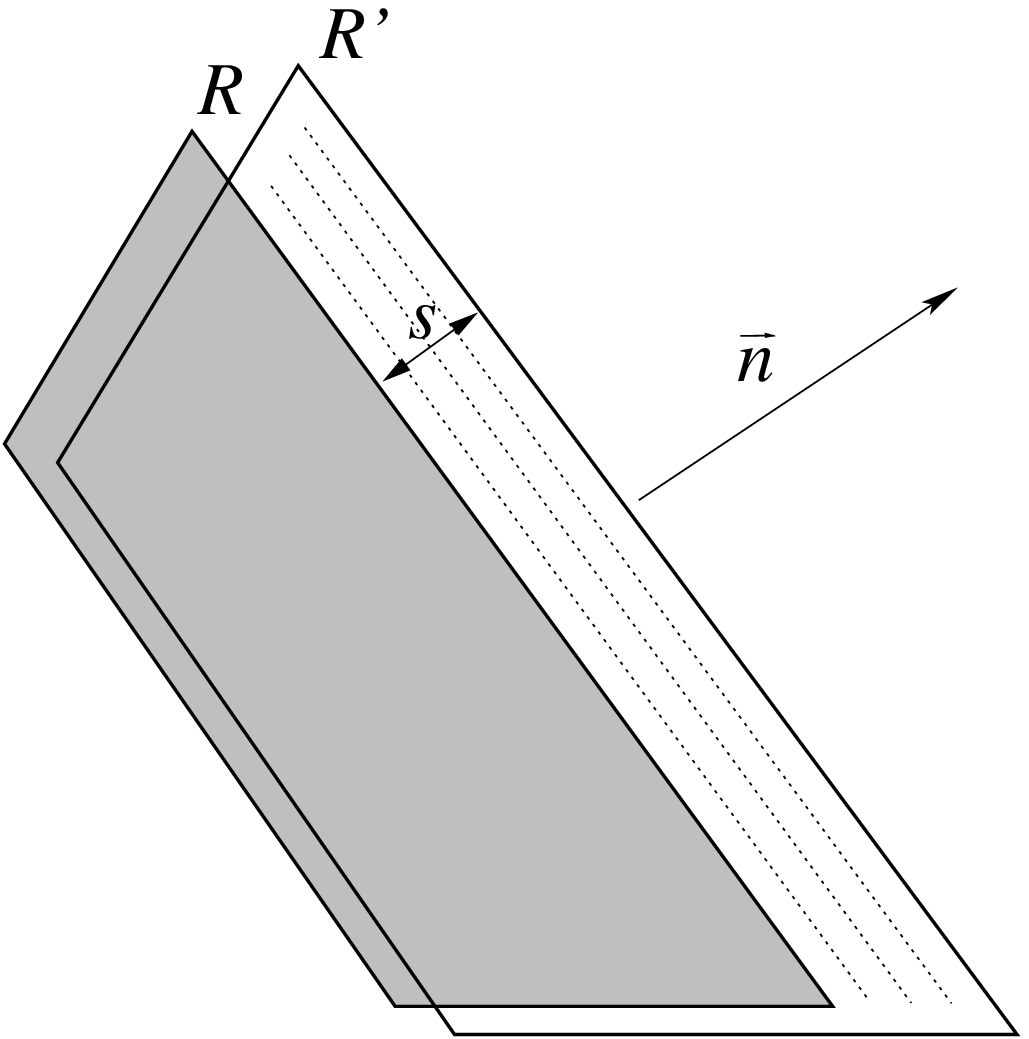}
\end{center}

First we find a hyperplane of one of the Kostant arrangements parallel to $F$ and outside the permutahedron. Recall from Remark~\ref{rem:outershift} that hyperplanes of the form~\eqref{eqn:layers} have nonnegative $\delta$-shifts. For positive $\delta$-shift, the corresponding hyperplane lies outside the permutahedron. In fact, we would like to maximize 
\begin{equation}
\mathrm{shift}(\{1,\ldots,j\}, W) = \sum_{i=1}^j\delta_{i} - \sum_{i=1}^j\delta_{w_i}
\end{equation}
because this will determine how much $p_R$ factorizes. The first sum is as large as possible because it is the sum of the first $j$ coordinates of $\delta=\frac{1}{2}(k-1,k-3,\ldots,-(k-3),-(k-1))$. Since $j\leq\lfloor k/2\rfloor$, we can pick $W$ disjoint from $\{1,\ldots,j\}$. Picking $W=\{k-j+1,\ldots,k\}$ means the second sum consists of the last (and smallest) entries of $\delta$. Thus $(k-1)/2$, $(k-3)/2$, $\ldots$, $(k-2j+1)/2$ appear in the first sum and their opposites in the second. The maximal shift is then 
\begin{equation}
\mathrm{shift}^{\mathrm{(max)}}(j) = 2\left(\frac{k-1}{2}+\frac{k-3}{2}+\cdots+\frac{k-2j+1}{2}\right) = j(k-j)\,.
\end{equation}
Suppose that $H(\lambda)$ is the hyperplane with this maximal shift (at distance $j(k-j)$ outside the permutahedron and parallel to $F$) and that it belongs to the Kostant arrangement $\mathcal{A}^{(\psi)}_\lambda$. 

The second step is to find a region $R'$ of $\mathcal{A}^{(\psi)}_\lambda$ with a facet on $H(\lambda)$ that overlaps with $R$. If we replace $\lambda$ by a multiple $m\lambda$ of itself, the partition of the permutahedron simply scales up by a factor of $m$, and the polynomials attached to the regions, as polynomials in $\lambda$ and $\beta$, remain the same (because the cells of the chamber complex $\mathcal{C}^{(k)}$ are cones). The hyperplanes of the Kostant arrangements almost scale, except for the $\delta$-shift factor. Those shifts preserve the distance between the hyperplanes and the ones supporting the facets of the permutahedron, even as the regions grow since the separation between parallel hyperplanes of $\mathcal{A}^{(\psi)}_\lambda$ increases. Hence for a large enough multiple of $\lambda$, one of the regions $R'$ of $\mathcal{A}^{(\psi)}_{m\lambda}$ with a facet on $H(m\lambda)$ will overlap with $mR$. From now on we'll assume that $\lambda$ has been replaced by a suitably large multiple of itself.

We are now in the setup of the above picture. The polynomials $p_R$ in $R$ and $P^{(\psi)}$ on $R'$ both give the multiplicities in the interior of their respective regions, and hence they are equal provided that $R\cap R'$ contains sufficiently many points. We can assume that we have scaled $\lambda$ sufficiently above so that this is the case. Since $P^{(\psi)}$ has to vanish outside of the permutahedron, it will vanish on the intersection of $R'$ with the hyperplanes
\begin{eqnarray*}
\beta_{\theta(1)}+\cdots+\beta_{\theta(j)} & = & \lambda_1+\cdots+\lambda_j + 1\\
 & \vdots &  \\
\beta_{\theta(1)}+\cdots+\beta_{\theta(j)} & = & \lambda_1+\cdots+\lambda_j + \mathrm{shift}^{\mathrm{(max)}}(j) - 1\,.
\end{eqnarray*}
If the intersection of $R'$ with these hyperplanes contains sufficiently many points (again, we can scale $\lambda$ so that this is the case), $P^{(\psi)}$ will have the defining equations of those hyperplanes as factors, and hence so will $p_R$. 

Here we have assumed $F$ is defined by \eqref{eqn:hyp1} and $\gamma=\gamma_{U,V}(\lambda)$ for $U=\theta(\{1,\ldots,j\})$, $V=\{1,\ldots,j\}$. If  $F$ is defined by \eqref{eqn:hyp2}, we get the same $U$ but $V=\{k-j+1,\ldots,k\}$, and the defining equations $\gamma, \gamma-1, \ldots, \gamma-\mathrm{shift}^{\mathrm{(max)}}(j)+1$. 
\end{proof}

We can lift this result to the weight polynomials associated to the cones of the chamber complex $\mathcal{C}^{(k)}$. This will allow us to think of the linear factors dividing the weight polynomials as polynomials both in $\lambda$ and $\beta$.

\begin{corollary}
\label{cor:lifting}
Let $\tau$ be the cone of $\mathcal{C}^{(k)}$ whose intersection with $L(\lambda)$ gives domain $R$ in the previous theorem, and $p_\tau$ its associated weight polynomial. If $\gamma_{U,V}(\lambda)+c$ divides $p_R$, then $\gamma_{U,V}+c$ divides $p_\tau$.
\end{corollary}

We will call these families of linear factors, \emph{parallel linear factors}. This shows that the smallest number of parallel linear factors is obtained when considering facets of the permutahedron normal to a permutation of $\omega_1$ or $\omega_{k-1}$.  In this case,  we get $k-2$ factors in the weight polynomials of the boundary regions on those facets. For $j=\lfloor k/2\rfloor$, we get a maximum of $\lfloor\frac{k}{2}\rfloor(k-\lfloor\frac{k}{2}\rfloor)-1 \sim k^2/4$ parallel linear factors. Since the $p_R$ have degree at most ${k-1\choose 2} \sim k^2/2$ in $\beta$ (regarding $\lambda$ as a parameter), we get linear factors accounting for about half the degree of the weight polynomials for those facets.

The fundamental weight $\omega_j = e_1+e_2+\cdots+e_j-\frac{j}{k}(1,1,\ldots,1)$ has an orbit of size ${k\choose j}$, and thus there are that many facets having a permutation of $\omega_j$ as an outer normal (the opposite parallel facets have the permutations of $\omega_{k-j}$ as normals). So the permutahedron for generic $\lambda = (\lambda_1,\ldots,\lambda_k)$ has $\sum_{j=1}^{k-1}{k\choose j} = 2^k-2$ facets, most of which have normals corresponding to central values of $j$ (i.e. close to $\lfloor k/2\rfloor$). The following table gives the minimum numbers of parallel linear factors for different values of $k$ and $j$. In parentheses are the numbers of facets having a permutation of $\omega_j$ or $\omega_{k-j}$ as a normal. For example, in $\sln{8}$, the maximal degree of the weight polynomials is 21 and we expect that the polynomials of regions with a facet on any of 112 of the 254 facets of the permutahedra to have 14 parallel linear factors.

\begin{displaymath}
\begin{array}{|l||c|c||c|c|c|c|}
\hline
& \#\textrm{(facets)} & \mathrm{max\,deg}(p_R) & j=1 & j=2 & j=3 & j=4 \\ 
\hline
\sln{3}\quad(A_2) & 6 & 1  & 1\ \hfill(6)  &         &          &           \\
\sln{4}\quad(A_3) & 14 & 3  & 2\ \hfill(8)  & 3 \ \hfill(6)   &          &           \\
\sln{5}\quad(A_4) & 30 & 6  & 3\ \hfill(10) & 5 \ \hfill(20)  &          &           \\
\sln{6}\quad(A_5) & 62 & 10 & 4\ \hfill(12) & 7 \ \hfill(30)  & 8 \ \hfill(20)   &          \\
\sln{7}\quad(A_6) & 126 & 15 & 5\ \hfill(14) & 9 \ \hfill(42)  & 11 \ \hfill(70)  &          \\
\sln{8}\quad(A_7) & 254 & 21 & 6\ \hfill(16) & 11 \ \hfill(56) & 14 \ \hfill(112) & 15 \ \hfill(70)  \\
\sln{9}\quad(A_8) & 510 & 28 & 7\ \hfill(18) & 13 \ \hfill(72) & 17 \ \hfill(168) & 19 \ \hfill(252) \\
\hline
\end{array}
\end{displaymath}

Theorem~\ref{thm:factorization} only depends on the fact that using our description in the previous section of the walls of the chamber complex for the Kostant partition function and the combinatorial description of the permutahedron (Lemma~\ref{lemma:KPFwalls}), we can argue that there will always be hyperplanes of Kostant arrangements parallel and close to the facets of the permutahedron. In order to extend the factorization phenomenon inside the permutahedron, we will need to use the complete description of the domains of polynomiality for the weight multiplicity function, obtained by symplectic geometry means in Theorem~\ref{thm:regular.pts}.

\subsection{Inside the permutahedron}

We already discussed at the end of the previous section that the hyperplanes supporting the walls partitioning the permutahedron are precisely the hyperplanes of the Kostant arrangements without the shift factors. We will take advantage here of overlaps between the improved domains of Section~\ref{sec:SymplectoStuff} and regions of the Kostant arrangements to show, not that the weight polynomials themselves factor, but rather that as we jump between two adjacent domains, the difference in the corresponding weight polynomials exhibits parallel linear factors. Given a facet between two adjacent domains of the permutahedron, we will see that we are able to find two hyperplanes of Kostant arrangements parallel to it and at maximal distance on either side of it, and deduce from this a number of parallel linear factors of the polynomial jump. 

\begin{definition}
We will say that two domains are \emph{adjacent} if they have the same dimension and a facet of one is a subset of a facet of the other, or equivalently if they intersect in a nonempty polytope of dimension one less.
\end{definition}

\begin{theorem}
\label{thm:InternalFactorization}
Let $P_1$ and $P_2$ be two adjacent full dimensional domains of polynomiality of the permutahedron for a generic dominant weight $\lambda$ of $\sln{k}$, and suppose that the normal to their touching facets is in the direction $\sigma(\omega_j)$ for some $\sigma\in\mathfrak{S}_k$. If $p_1$ and $p_2$ are the weight polynomials of $P_1$ and $P_2$, and $\gamma$ is the defining equation of the wall separating them, then the jump $p_1-p_2$ either vanishes or has the $j(k-j)-1$ linear factors 
\begin{displaymath}
(\gamma-s^{-}+1), (\gamma-s^{-}+2), \ldots, \gamma, \ldots, (\gamma+s^{+}-2), (\gamma+s^{+}-1)
\end{displaymath}
for some integers $s^{-}, s^{+} \geq 0$ satisfying
\begin{equation}
s^{-} + s^{+} = j(k-j)\,.
\end{equation}
\end{theorem}

\begin{proof}
Suppose the touching facets of $P_1$ and $P_2$ lie on the hyperplane 
\begin{displaymath}
\beta_{u_1}+\cdots+\beta_{u_j} = \lambda_{v_1}+\cdots+\lambda_{v_j}\,.
\end{displaymath}

Then among the Kostant arrangement hyperplanes 
\begin{displaymath}
\beta_{u_1}+\cdots+\beta_{u_j} = \lambda_{v_1}+\cdots+\lambda_{v_j} + \mathrm{shift}(V,W)\,,\qquad W\subseteq\{1,\ldots,k\}, \ |W|=j\,,
\end{displaymath}
we can find a pair for which the $\delta$-shift is minimal and maximal by picking appropriate subsets $W$. Clearly, the minimal shift $-s^{-}$ will be nonpositive, and the maximal shift $s^{+}$, nonnegative. In fact,
\begin{displaymath}
\begin{array}{rcccl}
s^{+} & = & \displaystyle\max_W\left(\sum_{i=1}^j\delta_{v_i}-\sum_{i=1}^j\delta_{w_i}\right) & = & \displaystyle\sum_{i=1}^j\delta_{v_i}-\min_{W}\sum_{i=1}^j\delta_{w_i}\\[5mm]
s^{-} & = & \displaystyle-\min_W\left(\sum_{i=1}^j\delta_{v_i}-\sum_{i=1}^j\delta_{w_i}\right) & = & \displaystyle\max_{W}\sum_{i=1}^j\delta_{w_i} - \sum_{i=1}^j\delta_{v_i}
\end{array}
\end{displaymath}
so that
\begin{equation}
s^{+}+s^{-} = \max_{W}\sum_{i=1}^j\delta_{w_i} - \min_{W}\sum_{i=1}^j\delta_{w_i} = 2\max_{W}\sum_{i=1}^j\delta_{w_i} = j(k-j)
\end{equation}
since $\delta=-\delta^{\mathrm{reverse}}$. For $k$ odd, the $\delta_i$ are integral, and hence so are $s^{-}$ and $s^{+}$. When $k$ is even, the $\delta_i$ are half-integers with odd numerators. Since we are adding/subtracting an even number of them ($2j$) to compute the shifts, we again get that $s^{-}$ and $s^{+}$ are integers.

We can find regions $Q_1$ and $Q_2$, $R_1$ and $R_2$ of Kostant arrangements as in the following diagram. We will think of these regions as open convex polytopes because in Lemma~\ref{lemma:KostantPolynomial} the polynomials giving the multiplicities on the regions of the Kostant arrangements are only valid in the interior of the regions. 
\begin{center}
\includegraphics[width=0.4\textwidth]{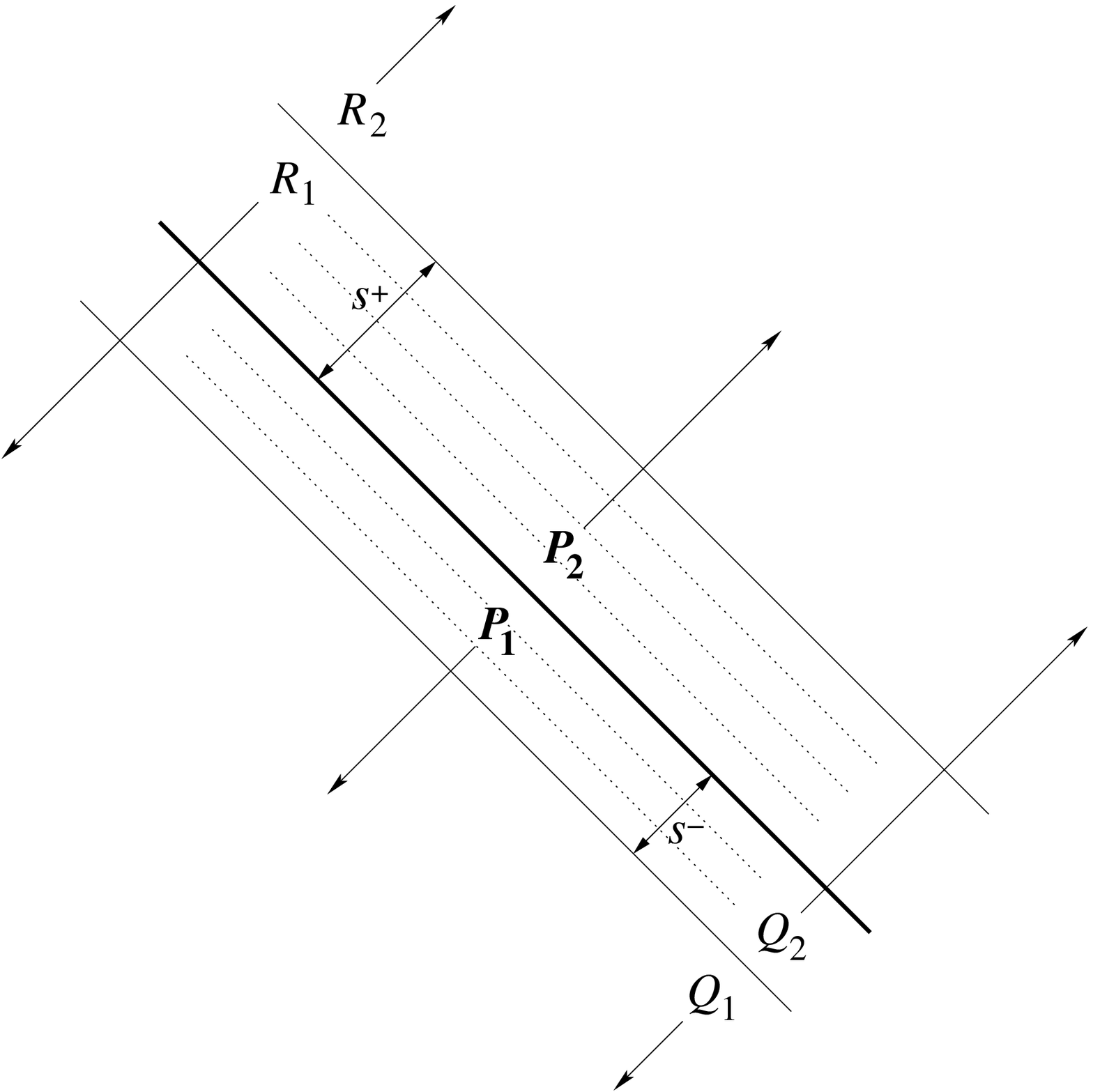}
\end{center}

We will let the corresponding polynomials, as given by $P^{(\psi)}$ in equation~(\ref{eqn:TypePolynomial}), be $q_1$ and $q_2$, $r_1$ and $r_2$ respectively. Since we can assume that we have scaled $\lambda$ sufficiently (as in Theorem~{\ref{thm:factorization}}), we have that $q_1=p_1=r_1$ and $q_2=p_2=r_2$, since $Q_1\cap P_1\cap R_1$ and $Q_2\cap P_2\cap R_2$ are large. Furthermore, $p_1$ and $q_2$ agree on $P_1\cap Q_2$ and similarly, $p_2$ and $r_1$ agree on $P_2\cap R_1$. Since $P_1\cap Q_2$ and $P_2\cap R_1$ contain enough lattice points on the bounded hyperplanes (dotted lines in the diagram), the differences $p_1-q_2$ and $p_2-r_1$ have to vanish on those hyperplanes. Hence  
\begin{equation} 
\begin{array}{rcl}
p_1-q_2 & = & (\gamma-s^{-}+1)(\gamma-s^{-}+2)\cdots(\gamma-1)\cdot h_1 \\[2mm]
p_2-r_1 & = & (\gamma+1)(\gamma+2)\cdots(\gamma+s^{+}-1)\cdot h_2 
\end{array}
\end{equation}
for some polynomials $h_1$ and $h_2$, unless $p_1=q_2$ or $p_2=r_1$, in which case $p_1=p_2$, since $p_1=r_1$ and $p_2=q_2$. If we assume that $p_1\neq p_2$, we have that
\begin{equation}
\begin{array}{rcl}
p_1-p_2 & = & (\gamma-s^{-}+1)(\gamma_s^{-}+2)\cdots(\gamma-1)\cdot h_1 \\[2mm]
p_2-p_1 & = & (\gamma+1)(\gamma+2)\cdots(\gamma+s^{+}-1)\cdot h_2 
\end{array}
\end{equation}
and since $p_1$ and $p_2$ have to agree on the lattice points on the wall between $P_1$ and $P_2$, their difference is also divisible by $\gamma$. Hence we get
\begin{equation}
p_1-p_2 = (\gamma-s^{-})(\gamma-s^{-}+1)\cdots\gamma\cdots(\gamma+s^{+}-1)(\gamma+s^{+})\cdot h_3
\end{equation}
for some $h_3$.
\end{proof}

\begin{remark}
As in Corollary~\ref{cor:lifting}, we can lift this result to the weight polynomials $p_\tau$ associated to the cones of the chamber complex $\mathcal{C}^{(k)}$ and regard the parallel linear factors as polynomials in both $\lambda$ and $\beta$.
\end{remark}

\section{The chamber complexes for $\sln{3}$ ($A_2$) and $\sln{4}$ ($A_3$)}
\label{sec:A2A3}

In this section we explicitly compute the chamber complexes for $k=3$ and $k=4$.  For $k=4$, we find that the chamber complex does not optimally partition the domains of polynomiality for the multiplicity function.  In Theorem~\ref{thm:A3partition}, we prove that the optimal glued complex does agree with Theorem~\ref{thm:SamePartition} for $k=4$.

\subsection{The chamber complex for $\sln{3}$ ($A_2$)}

Using the procedure described in Section~\ref{sec:SPF} to write down the multiplicity function as a single partition function in the case $k=3$ ($A_2$) gives that
\begin{equation}
m_\lambda(\beta) = \phi_{E_3}\left(B_3{\lambda\choose\beta}\right)
\end{equation}
with
\begin{equation}
\renewcommand{\arraystretch}{1.0}
E_3 = \left(\begin{array}{rrrrrr}
1 &  1 & 0 & 0 & 0 & 0 \\
0 & -1 & 1 & 0 & 0 & 0 \\
0 &  1 & 0 & 1 & 0 & 0 \\
0 & -1 & 0 & 0 & 1 & 0 \\
0 & -1 & 0 & 0 & 0 & 1
\end{array}\right) \quad \textrm{and} \quad 
B_3{\lambda\choose\beta} = \left(\begin{array}{c}
\lambda_1 - \lambda_2\\ 2\lambda_2 - \beta_1 - \beta_2\\ \beta_1 + \beta_2 + \lambda_1\\ \lambda_2 - \beta_1\\ \lambda_2 - \beta_2\end{array}\right)\,. 
\end{equation}

We can compute the chamber complex associated to $E_3$ and intersect it with the space 
\begin{equation}
\tilde{B} = \left\{B_3{\lambda\choose\beta}\ : \ \lambda\in\mathbb{R}^3, \beta\in\mathbb{R}^3, \lambda_1+\lambda_2+\lambda_3=0, \beta_1+\beta_2+\beta_3=0\right\}\,.
\end{equation}

In that space we can apply
\begin{equation}
B_3^{-1} = \frac{1}{9}\left(\begin{array}{rrrrr}
6 & 2 & 3 & 1 & 1\\
-3 & 2 & 3 & 1 & 1\\
-3 & -1 & 3 & -5 & 4\\
-3 &-1 &3 &4 &-5
\end{array}\right)
\end{equation}
to rectify the cones of that complex to obtain $\mathcal{C}^{(3)}$. The full dimensional cones of $\mathcal{C}^{(3)}$ are given by
\begin{equation}
\renewcommand{\arraystretch}{1.8}
\begin{array}{|l|l|l|l|}
\hline
& & \tau_3 = \mathrm{pos}(b,a_1,c_2,c_3) & \tau_6 = \mathrm{pos}(b,a_1,a_2,c_3) \\
\tau_1 = \mathrm{pos}(b,a_1,a_2,a_3) & \tau_2 = \mathrm{pos}(b,c_1,c_2,c_3) & \tau_4 = \mathrm{pos}(b,a_2,c_1,c_3) & \tau_7 = \mathrm{pos}(b,a_1,a_3,c_2) \\
& & \tau_5 = \mathrm{pos}(b,a_3,c_1,c_2) & \tau_8 = \mathrm{pos}(b,a_2,a_3,c_1) \\
\hline
\end{array}
\end{equation}
where the rays are
\begin{equation}
\begin{array}{l@{\hspace{10mm}}l@{\hspace{10mm}}l}
a_1 = [2,-1,-1,2,-1,-1] &                      & c_1 = [1,1,-2,-2,1,1]\\
a_2 = [2,-1,-1,-1,2,-1] & b   = [1,0,-1,0,0,0] & c_2 = [1,1,-2,1,-2,1]\\
a_3 = [2,-1,-1,-1,-1,2] &                      & c_3 = [1,1,-2,1,1,-2]\\
\end{array}
\end{equation}

The cones above are grouped into orbits under the action of the symmetric group $\mathfrak{S}_3$ on the $\beta$-coordinates. In general the set of cones of $\mathcal{C}^{(k)}$ won't be closed under the action of $\mathfrak{S}_k$ on the $\beta$-coordinates, even though the multiplicities under the $\mathfrak{S}_k$ should be invariant. 

We can get the polynomial $p_i$ corresponding to $\tau_i$ easily through interpolation, using for example the Kostant partition function for $A_2$ which has the simple form
\begin{displaymath}
\renewcommand{\arraystretch}{1.0}
K(a,-a+b,-b) = \left\{\begin{array}{ll}
\mathrm{min}\{a,b\}+1 & \quad\textrm{if $a,b\in\mathbb{N}$}\,,\\
0 & \quad\textrm{otherwise.}
\end{array}\right.
\end{displaymath}

\begin{displaymath}
\renewcommand{\arraystretch}{1.8}
\begin{array}{|l|l|l|l|}
\hline
& & p_3 = 1+\lambda_1-\beta_1 & p_6 = 1+\lambda_1+\lambda_2+\beta_3 \\
p_1 = 1+\lambda_2-\lambda_3 & p_2 = 1+\lambda_1-\lambda_2 & p_4 = 1+\lambda_1-\beta_2 & p_7 = 1+\lambda_1+\lambda_2+\beta_2 \\
& & p_5 = 1+\lambda_1-\beta_3 & p_8 = 1+\lambda_1+\lambda_2+\beta_1 \\
\hline
\end{array}
\end{displaymath}

Note that even though they highlight the symmetries in the $\beta_i$'s, these polynomials are a little ambiguous since they are defined up to the relations $\lambda_1+\lambda_2+\lambda_3=0$ and $\beta_1+\beta_2+\beta_3=0$, which allow for some substitutions to be made. To avoid any ambiguity, we can rewrite them in terms of the fundamental weight basis $\omega_1=\frac{1}{3}(2,-1,-1)$ and $\omega_2=\frac{1}{3}(1,1,-2)$. Then if $\lambda=l_1\omega_1+l_2\omega_2$ and $\beta=b_1\omega_1+b_2\omega_2$, the polynomials take the form
\begin{displaymath}
\renewcommand{\arraystretch}{1.8}
\begin{array}{|l|l|l|l|}
\hline
& & p_3 = 1+\frac{1}{3}(2l_1+l_2-2b_1-b_2) & p_6 = 1+\frac{1}{3}(l_1+2l_2-b_1-2b_2) \\
p_1 = 1+l_2 & p_2 = 1+l_1 & p_4 = 1+\frac{1}{3}(2l_1+l_2+b_1-b_2) & p_7 = 1+\frac{1}{3}(l_1+2l_2-b_1+b_2) \\
& & p_5 = 1+\frac{1}{3}(2l_1+l_2+b_1+2b_2) & p_8 = 1+\frac{1}{3}(l_1+2l_2+2b_1+b_2) \\
\hline
\end{array}
\end{displaymath}

\begin{figure}[!ht]
\begin{center}
\includegraphics[width=1.0\textwidth]{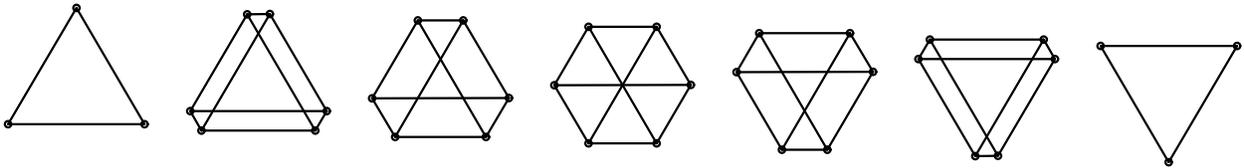}
\caption{Weight diagrams and their domains of polynomiality for $A_2$.}
\label{fig:A2diagrams}
\end{center}
\end{figure}

The domains of polynomiality of a weight diagram for a given $\lambda$ will be the (possibly empty) polytopes $\tau_i\cap L(\lambda)$ for $i=1,2,\ldots,8$, so there are at most eight of them, although in practice at most seven appear at a time. We could obtain a symbolic description of the domains of polynomiality for the weight diagram of any $\lambda$ of $A_2$ from the chamber complex. When $\lambda$ is one of the fundamental weights, we get a triangle with constant multiplicities inside; otherwise we get a hexagon with a (possibly empty) central triangle in which the multiplicities are constant and decrease linearly outside. Figure~\ref{fig:A2diagrams} shows what happens when we move from one fundamental weight to the other. The picture for $A_2$ is already well-known (see, for example, \cite{Humphreys,BG}). 

Corollary~\ref{cor:LambdaComplex} explains why the second and third diagrams, as well as the fifth and sixth of Figure~\ref{fig:A2diagrams} are variations of each other, and why the polynomials attached to each of the seven regions are the same for each of these pairs of diagrams. For $A_2$, the cones $\tau_k$ project under $p_\Lambda$ to the three cones
\begin{displaymath}
\begin{array}{lclcl}
C_0 & = & \mathrm{pos}((2,-1,-1),(1,1,-2)) & = & \mathrm{pos}(\omega_1,\omega_2)\\
C_1 & = & \mathrm{pos}((2,-1,-1),(1,0,-1)) & = & \mathrm{pos}(\omega_1,\omega_1+\omega_2)\\
C_2 & = & \mathrm{pos}((1,0,-1),(1,1,-2)) & = & \mathrm{pos}(\omega_1+\omega_2,\omega_2)\,.
\end{array}
\end{displaymath}

We can see that $C_1$ and $C_2$ partition the fundamental Weyl chamber $C_0$ of $A_2$, and hence $\mathcal{C}^{(2)}_\Lambda$ consists of $C_1$, $C_2$ and all their faces. Therefore for $A_2$ there are only two generic types of $\lambda$'s : $\lambda$'s with $\lambda_2<0$ (diagrams 2 and 3 on Figure~\ref{fig:A2diagrams}) and $\lambda$'s with $\lambda_2>0$ (diagrams 5 and 6 on Figure~\ref{fig:A2diagrams}). The case $\lambda_2=0$ corresponds to the regular hexagon, while the degenerate cases $\lambda_1=\lambda_2$ and $\lambda_2=\lambda_3$ correspond to the triangles. If we express $\lambda$ in terms of the fundamental weights $\lambda=l_1\omega_1+l_2\omega_2$, these correspond to $l_1<l_2$, $l_1>l_2$ and $l_1=l_2$ respectively for the hexagons, and $l_1=0$ and $l_2=0$ for the degenerate cases (triangles). 

\subsection{The chamber complex for $\sln{4}$ ($A_3$)}

We can write
\begin{displaymath}
m_\lambda(\beta) = \phi_{E_4}\left(B_4{\lambda\choose\beta}\right)
\end{displaymath}
with
\begin{displaymath}
\renewcommand{\arraystretch}{1.2}
E_4 = \left(\begin{array}{rrrrrrrrrrrr}
0& 1& 1& 1& 0& 0& 0& 0& 0& 0& 0& 0\\
1& -1& 1& 0& 1& 0& 0& 0& 0& 0& 0& 0\\ 
0& 1& 0& 0& 0& 1& 0& 0& 0& 0& 0& 0\\ 
0& 0& 1& 0& 0& 0& 1& 0& 0& 0& 0& 0\\ 
0& -1& -1& 0& 0& 0& 0& 1& 0& 0& 0& 0\\ 
-1& 0& -2& 0& 0& 0& 0& 0& 1& 0& 0& 0\\ 
-1& 0& -1& 0& 0& 0& 0& 0& 0& 1& 0& 0\\ 
-1& 0& -1& 0& 0& 0& 0& 0& 0& 0& 1& 0\\ 
1& -1& 0& 0& 0& 0& 0& 0& 0& 0& 0& 1
\end{array}\right) \quad \textrm{and} \quad B_4{\lambda\choose\beta} = 
\left(\begin{array}{c}
\lambda_1 +\beta_1 +\beta_2 +\beta_3 \\ \lambda_2 -\lambda_3 \\ \lambda_1 -\lambda_2 \\ \lambda_2 -\lambda_3 \\ \lambda_2 +2\lambda_3-\beta_1 -\beta_2 -\beta_3 \\ 2\lambda_3-\beta_1 -\beta_2 \\ \lambda_3 -\beta_1 \\ \lambda_3 -\beta_2 \\ \lambda_2 -\beta_3 
\end{array}\right)\,.
\end{displaymath}

\begin{remark}
$E_4$ is not unimodular. We don't know of a unimodular matrix for $\sln{4}$ that would make the multiplicity function into a single partition function.
\end{remark}

For $A_3$ we must use the computer to do most of the computations. A symbolic calculator like Maple or Mathematica is especially useful. Here we used Maple (versions 7 and 8) and the package $\texttt{convex}$ by Matthias Franz \cite{convex}.

The set $\mathcal{B}^{(4)}$ of bases for $E_4$ has 146 elements, so there are 146 base cones $\tau_\sigma$ for $\sigma\in\mathcal{B}^{(4)}$. These are 9-dimensional cones, however they collapse to 132 6-dimensional cones when intersected with 
\begin{displaymath}
\tilde{B} = \left\{B_4{\lambda\choose\beta}\ : \ \lambda\in\mathbb{R}^4, \beta\in\mathbb{R}^4, \lambda_1+\cdots+\lambda_4=0, \beta_1+\cdots+\beta_4=0\right\}
\end{displaymath}

The full chamber complex is the complex of all intersections of the 146 cones in $\mathbb{R}^9$, and it has 6472 full-dimensional cones. The chamber complex in $\tilde{B}$-space has 1202 full-dimensional (6-dimensional) cones $\tilde{\tau}_k$, which we can rectify to get the chamber complex $\mathcal{C}^{(4)}$ in $(\lambda,\beta)$-space with cones $\tau_k$, $k=1, \ldots, 1202$. However, the 6-dimensional chamber complex thus obtained is not closed under the action of the symmetric group $\mathfrak{S}_4$ on the $\beta$-coordinates.

Despite the fact that the chamber complex seems to lack the symmetry property in $\beta$, we will see, as we find the polynomials attached to the domains of polynomiality, that there is a way to regain it. We can compute the polynomial associated to each of these 1202 cones by interpolation, for example using the fact that De~Loera and Sturmfels computed the polynomials for the Kostant partition function for $A_3$ in \cite{DeLoeraSturmfels}. These 1202 polynomials are not all distinct.

\begin{observation}
\label{rem:GluedComplex}
If we group together the top-dimensional cones from $\{\tau_k\ : \ k=1,\,\ldots,\,1202\}$ with a particular polynomial, their union is always a convex polyhedral cone again. Grouping cones this way yields a glued chamber complex $\mathcal{G}$ in $(\lambda,\beta)$-space with 612 cones $G_k$, $k=1,\,\ldots,\,612$. These cones form 64 orbits under the action of $\mathfrak{S}_4$ on the $\beta$-coordinates.
\end{observation}

\begin{proof}
Here is a description of the algorithm used to make this observation. Suppose that $\{\tau_{i_1}, \tau_{i_2}, \ldots, \tau_{i_N}\}$ consists of all the cones with a particular given associated polynomial, and let $\tau$ be the convex polyhedral cone spanned by the union of all their rays. We want to prove $\tau=\cup_{j=1}^N\tau_{i_j}$.

We can find an affine half-space whose intersection with each of these cones is non-empty and bounded, so that we can work with truncated cones. The half-space $\lambda_1\leq1$ works.  The union of $\{\tau_{i_1}, \ldots, \tau_{i_N}\}$ will equal $\tau$ if and only if the union of their truncations gives the truncation of $\tau$. The truncated cones are polytopes, and we can compute their volume. We can check that the union of the truncations of $\tau_{i_1}, \ldots, \tau_{i_N}$ is the truncation of $\tau$ just by checking that the volumes match. We know the $\tau_{i_j}$ have disjoint interiors because they are defined as the common refinement of base cones, hence the volume of the union of all these truncated cones is simply the sum of their volumes. If the computations are done symbolically (in Maple), there is no danger that truncated cones with very small volumes could create round-off errors.

The volumes are compared symbolically for every family of cones corresponding to the same polynomial. We glue together all the cones with the same polynomial, and observe that we still get a complex of convex polyhedral cones. This glued complex is invariant under the action of $\mathcal{S}_4$ on the $\beta$-coordinates. 
\end{proof}

We now have two 6-dimensional chamber complexes $\mathcal{C}^{(3)}$ and $\mathcal{G}$, and we can construct the complexes $\mathcal{C}^{(3)}_\Lambda$ and $\mathcal{G}_\Lambda$ by first projecting all the cones through $p_\Lambda$ and then forming their common refinement.

After transporting the hyperplane $x+y+z+w=0$ of $\mathbb{R}^4$ into the hyperplane $z=0$ through the orientation-preserving isometry
\begin{equation}
\renewcommand{\arraystretch}{1.0}
T_4 := \frac{1}{2}\left(\begin{array}{rrrr}
1 & -1 & -1 & 1 \\
1 & -1 & 1 & -1 \\
1 & 1 & -1 & -1 \\
-1 & -1 & -1 & -1
\end{array}\right)\,,
\end{equation}
we can work in the coordinates $(x,y,z)$ and look at the intersections of the complexes with the hyperplane $z=1$ of $\mathbb{R}^3$. 

\begin{figure}[!ht]
\begin{center}
\includegraphics[width=0.9\textwidth]{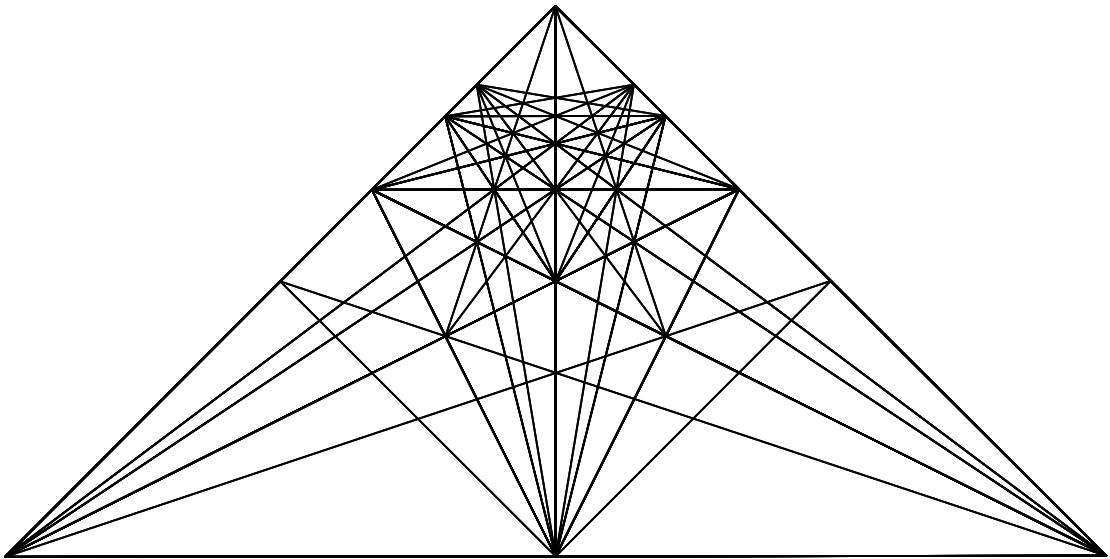}
\caption{The chamber complex $\mathcal{C}^{(3)}_\Lambda$.}
\label{fig:lambda_complex}
\end{center}
\end{figure}

\begin{figure}[!ht]
\begin{center}
\includegraphics[width=0.9\textwidth]{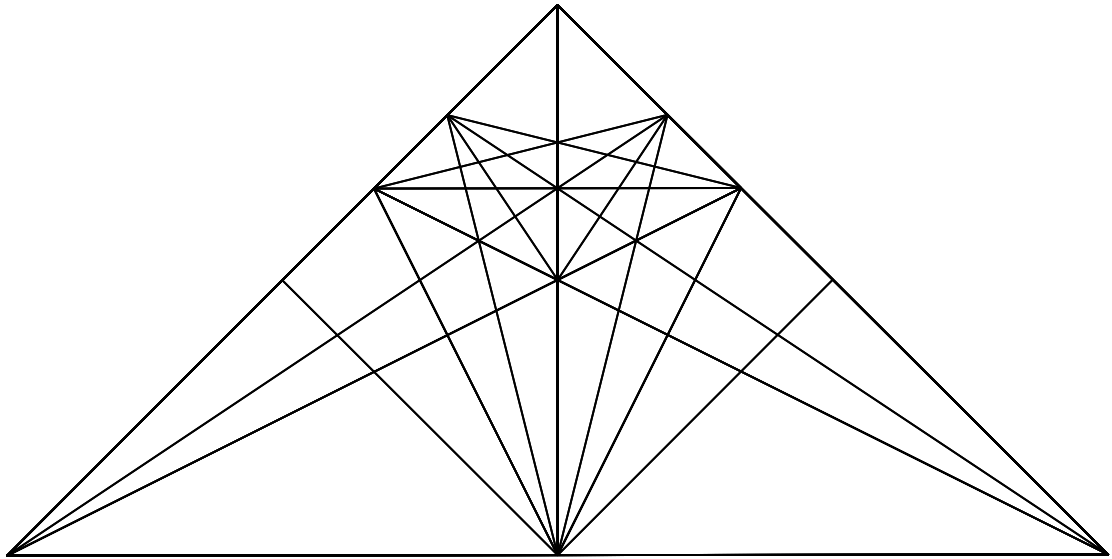}
\caption{The chamber complex $\mathcal{G}_\Lambda$.}
\label{fig:glued_lambda_complex}
\end{center}
\end{figure}

\begin{figure}[!ht]
\begin{center}
\includegraphics[width=0.6\textwidth]{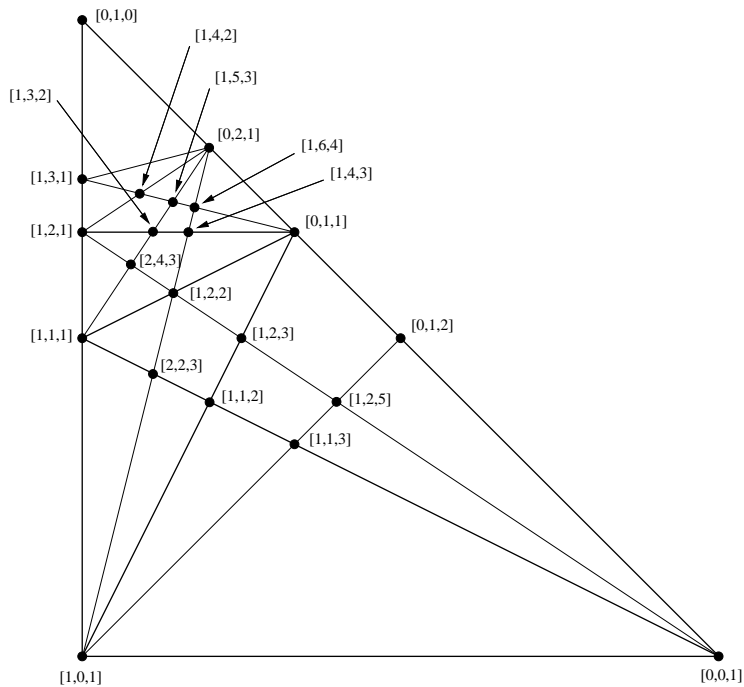}
\caption{The chamber complex $\mathcal{G}_\Lambda$ for $\lambda_1<-\lambda_4$ in terms of fundamental weights.}
\label{fig:half_glued_lambda_complex}
\end{center}
\end{figure}

\begin{figure}[!ht]
\begin{center}
\fboxsep50pt\fboxrule1pt\fbox{\includegraphics[width=0.4\textwidth]{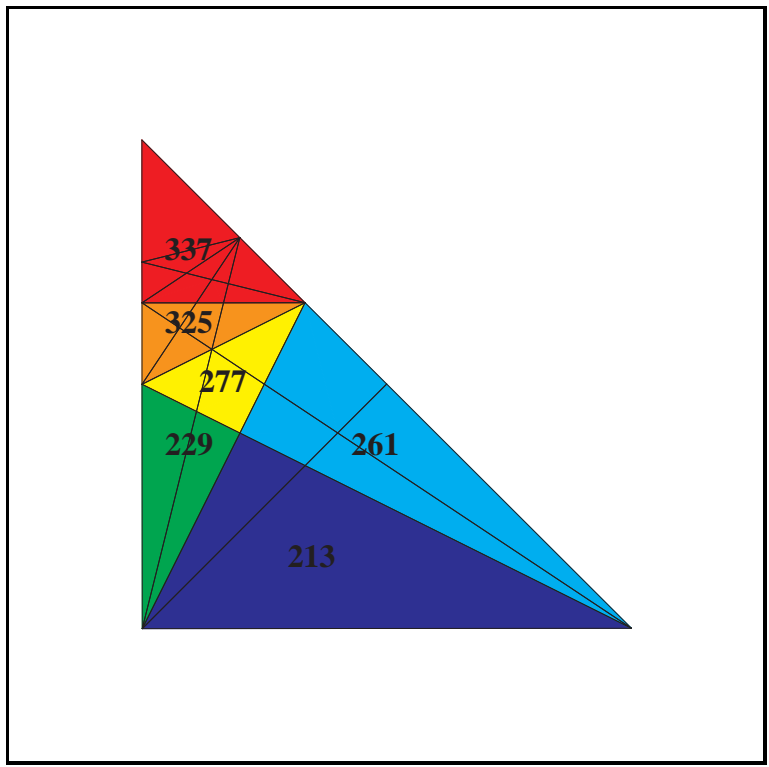}}
\caption{The numbers of regions in the permutahedra.}
\label{fig:number_regions_permutahedra}
\end{center}
\end{figure}

Figures~\ref{fig:lambda_complex}--\ref{fig:half_glued_lambda_complex} show the complexes intersected with the hyperplane $z=1$ and also the complexes modulo the symmetry $\lambda\mapsto-\lambda^{\mathrm{rev}}$, which translates into a reflection along the central (vertical) line of the complexes. This symmetry reflects the symmetry of the Dynkin diagrams for $A_n$. Figure~\ref{fig:number_regions_permutahedra} shows that even though regions appear and disappear along the lines of the complex (facets of the full-dimensional cones of the complexes), the complex given by simply looking at the number of regions in the permutahedra is coarser. 

\begin{observation}
For $A_3$, only six generic cases occur. Generic permutahedra are always partitioned into 213, 229, 261, 277, 325 or 337 regions. Degenerate cases occur along the walls in Figure~\ref{fig:glued_lambda_complex}. 
\end{observation}

Projecting the cones of the glued complex $\mathcal{G}$ on $\lambda$-space gives 62 distinct cones, 60 of them corresponding to individual orbits under the action of $\mathfrak{G}_4$ on the $\beta$-coordinates. The chamber complex $\mathcal{G}_\Lambda$ we get by taking their common refinement has 50 regions, or 25 modulo the symmetry $\lambda\mapsto -\lambda^{\mathrm{rev}}$. This complex classifies the combinatorial types of $\lambda$, i.e. the $\lambda$'s with the same partitioned permutahedra and family of polynomials.

We now give a simple proof of Theorem~\ref{thm:SamePartition} for $A_3$ and show that regions of the permutahedron given by \eqref{eqn:walls} are as large as possible.   

\begin{theorem}
\label{thm:A3partition}
For $A_3$, the optimal partition of the permutahedron into domains of polynomiality for the weight multiplicities coincides with the partition of the permutahedron into domains of polynomiality for the Duistermaat-Heckman measure.
\end{theorem}

\begin{proof}
We give a computer verified proof that for all $\lambda$ in the fundamental Weyl chamber of $\sln{4}$ the intersection of $\mathcal{G}$ with $L(\lambda)$ defines walls within the permutahedron as determined by Theorem~\ref{thm:regular.pts}. We do this by expressing the walls of the permutahedron as the convex hulls of subsets of its vertices. The following is an outline of our algorithm.

The full dimensional cones $G_1,\,\ldots,\,G_{612}$ of the complex $\mathcal{G}$, when intersected $L(\lambda)$, subdivide the permutahedron into regions. For generic $\lambda$, $G_k\cap L(\lambda)$ is either empty or a 3-dimensional region of the permutahedron. Furthermore, a 2-dimensional facet of that region will come from the intersection of a facet $F$ of $G_k$ with $L(\lambda)$, and an edge of that 2-dimensional facet will come from the intersection of a facet $L$ of $F$ with $L(\lambda)$, and finally a vertex of that edge will come from the intersection of a facet of $L$ with $L(\lambda)$. 

\begin{enumerate}
\item Set $\mathcal{F}$ equal to the set of all facets of the cones $G_1,\,\ldots,\,G_{612}$.  

\item Classify the facets in $\mathcal{F}$ according to their normals: call $\mathcal{F}_i$ the subset of $\mathcal{F}$ consisting of all the facets with normal direction $n_i$. Since each facet lies on a unique hyperplane, and since all these hyperplanes go through the origin, two facets will lie on the same hyperplane if and only if they have the same normals up to a scalar multiple. In our case, we find that there are 37 distinct normal directions. 

\item Set $K_i=\cup_{F\in\mathcal{F}_i}F$ and verify that $K_i$ is again a convex polyhedral cone. The verification is done by a truncation and volume comparison method similar to the one used in Observation~\ref{rem:GluedComplex}.  The intersection of the $K_i$ with $L(\lambda)$ will be the walls partitioning the permutahedron. 

\item For each $i$, set $\mathcal{V}_i$ to be the set of facets of facets of $K_i$. The elements of $\mathcal{V}_i$ are three dimensional cones. 

\item For each $i$, identify the $f\in\mathcal{V}_i$ whose intersection with $L(\lambda)$ for generic $\lambda$ is a point. The convex hull of those points is $K_i\cap L(\lambda)$. These points are all vertices of the permutahedron, and the walls they define are exactly those of equation~\eqref{eqn:walls}.
\end{enumerate}

We will illustrate this last step on an example. We find that one of the $\mathcal{V}_i$ consists of the 10 cones, which we will denote $f_{1},\,\ldots,\,f_{10}$. One remarkable thing about the cones $f_j$ is that the first four coordinates of their rays always correspond to one of the fundamental weights, while the last four correspond to a conjugate of the same fundamental weight. That is true for all $\mathcal{V}_i$. We have
\begin{displaymath}
\renewcommand{\arraystretch}{1.5}
\begin{array}{lcllcl}
f_1 & = & \mathrm{pos}(\omega_1, \omega_1),\  (\omega_2, \pi\cdot\omega_2),\  (\omega_3, \omega_3)) & f_6 & = & \mathrm{pos}((\omega_2, \pi\cdot\omega_2),\  (\omega_3, \phi\cdot\omega_3),\  (\omega_3, \omega_3))\\
f_2 & = & \mathrm{pos}((\omega_1, \omega_1),\  (\omega_2, \pi\cdot\omega_2),\  (\omega_3, \phi\cdot\omega_3)) & f_7 & = & \mathrm{pos}((\omega_1, \sigma\cdot\omega_1),\  (\omega_1, \omega_1),\  (\omega_3, \phi\cdot\omega_3))\\
f_3 & = & \mathrm{pos}((\omega_1, \sigma\cdot\omega_1),\  (\omega_2, \pi\cdot\omega_2),\  (\omega_3, \phi\cdot\omega_3))& f_8 & = & \mathrm{pos}((\omega_1, \omega_1),\  (\omega_3, \phi\cdot\omega_3),\  (\omega_3, \omega_3))\\
f_4 & = & \mathrm{pos}((\omega_1, \sigma\cdot\omega_1),\  (\omega_2, \pi\cdot\omega_2),\  (\omega_3, \omega_3) & f_9 & = & \mathrm{pos}((\omega_1, \sigma\cdot\omega_1),\  (\omega_1, \omega_1),\  (\omega_3, \omega_3))\\
f_5 & = & \mathrm{pos}((\omega_1, \sigma\cdot\omega_1),\  (\omega_3, \phi\cdot\omega_3),\  (\omega_3, \omega_3)) & f_{10}&=& \mathrm{pos}((\omega_1, \sigma\cdot\omega_1),\  (\omega_1, \omega_1),\  (\omega_2, \pi\cdot\omega_2))
\end{array}
\end{displaymath}
where $\sigma=(1\ 3)$, $\pi=(2\ 3)$, $\phi=(2\ 4)$.

To find the intersection of one of these cones with $L(\lambda)$, we want to see whether there is a linear combination of its rays with nonnegative coefficients that would lie in $L(\lambda)$. If the rays are $r_1,\,\ldots,\,r_s$, we are looking for $a_1,\,\ldots,\,a_s\geq0$ such that
\begin{displaymath}
a_1\,r_1 + a_2\,r_2 + \cdots + a_s\,r_s = (\lambda_1,\lambda_2,\lambda_3,\lambda_4,*,*,*,*)\,,
\end{displaymath}
or equivalently,
\begin{displaymath}
a_1\,p_\Lambda(r_1) + a_2\,p_\Lambda(r_2) + \cdots + p_\Lambda(a_s\,r_s) = (\lambda_1,\lambda_2,\lambda_3,\lambda_4)\,,
\end{displaymath}

Hence we will get vertices for those $\lambda$'s and $f_j$'s such that $\lambda\in p_\Lambda(f_j)$. So we compute the $p_\Lambda(f_j)$:
\begin{displaymath}
\renewcommand{\arraystretch}{1.2}
\begin{array}{lcl@{\hspace{0.1\textwidth}}lcl}
p_\Lambda(f_1) & = & \mathrm{pos}(\omega_1, \omega_2, \omega_3) & p_\Lambda(f_5) & = & \mathrm{pos}(\omega_1, \omega_3)\\
p_\Lambda(f_2) & = & \mathrm{pos}(\omega_1, \omega_2, \omega_3) & p_\Lambda(f_6) & = & \mathrm{pos}(\omega_2, \omega_3)\\
p_\Lambda(f_3) & = & \mathrm{pos}(\omega_1, \omega_2, \omega_3) & p_\Lambda(f_7) & = & \mathrm{pos}(\omega_1, \omega_3)\\
p_\Lambda(f_4) & = & \mathrm{pos}(\omega_1, \omega_2, \omega_3) & p_\Lambda(f_8) & = & \mathrm{pos}(\omega_1, \omega_3)\\
&&& p_\Lambda(f_9) & = & \mathrm{pos}(\omega_1, \omega_3)\\
&&& p_\Lambda(f_{10})&=& \mathrm{pos}(\omega_1, \omega_2)
\end{array}
\end{displaymath}

Only the first four of the cones span the fundamental Weyl chamber; the other six won't intersect $L(\lambda)$ for generic $\lambda$. 

Observing that the last four coordinates of the rays of the $f_j$'s can always be obtained by applying a single permutation to the first four, we can rewrite $f_1,\,f_2,\,f_3,\,f_4$ as
\begin{eqnarray*}
f_1 & = & \mathrm{pos}((\omega_1,(2\,3)\cdot\omega_1),\ (\omega_2,(2\,3)\cdot\omega_2),\ (\omega_3,(2\,3)\cdot\omega_3)) \\
f_2 & = & \mathrm{pos}((\omega_1,(2\,4\,3)\cdot\omega_1),\ (\omega_2,(2\,4\,3)\cdot\omega_2),\ (\omega_3,(2\,4\,3)\cdot\omega_3)) \\
f_3 & = & \mathrm{pos}((\omega_1,(1\,2\,4\,3)\cdot\omega_1),\ (\omega_2,(1\,2\,4\,3)\cdot\omega_2),\ (\omega_3,(1\,2\,4\,3)\cdot\omega_3)) \\
f_4 & = & \mathrm{pos}((\omega_1,(1\,2\,3)\cdot\omega_1),\ (\omega_2,(1\,2\,3)\cdot\omega_2),\  (\omega_3,(1\,2\,3)\cdot\omega_3))
\end{eqnarray*}

It then follows that
\begin{eqnarray*}
f_1\cap L(\lambda) & = & (\lambda, (2\,3)\cdot\lambda) \\
f_2\cap L(\lambda) & = & (\lambda, (2\,4\,3)\cdot\lambda) \\
f_3\cap L(\lambda) & = & (\lambda, (1\,2\,4\,3)\cdot\lambda) \\
f_4\cap L(\lambda) & = & (\lambda, (1\,2\,3)\cdot\lambda)
\end{eqnarray*} 
which means there will be a wall with vertices
\begin{displaymath}
\begin{array}{ccccr}
(2\,3)\cdot\lambda       & = & (\lambda_1,\lambda_3,\lambda_2,\lambda_4) & = & \lambda'\\
(1\,2\,3)\cdot\lambda    & = & (\lambda_3,\lambda_1,\lambda_2,\lambda_4) & = & (1\,2)\lambda'\\
(2\,4\,3)\cdot\lambda    & = & (\lambda_1,\lambda_3,\lambda_4,\lambda_2) & = & (3\,4)\lambda'\\
(1\,2\,4\,3)\cdot\lambda & = & (\lambda_3,\lambda_1,\lambda_4,\lambda_2) & = & (1\,2)(3\,4)\lambda'
\end{array}
\end{displaymath}
in the permutahedron for $\lambda$. This wall is the convex hull of $W\cdot\lambda'$ with $W$ the parabolic subgroup $\langle (1\ 2), (3\ 4)\rangle$. All these vertices have the same scalar product with $(1,1,-1,-1)$, and thus they lie on the same hyperplane with normal $(1,1,-1,-1)$. They are also the only points of $\mathcal{S}_4\cdot\lambda$ lying on that hyperplane (for $\lambda$ generic).

This process is automated, and we can express all the walls in this fashion by repeating this process for all the $\mathcal{V}_i$'s. We finally check that these walls are the same as those that partition the permutahedron for the DH-measure.
\end{proof}

The set of domains for a given permutahedron is closed under the action of the Weyl group, so they come into orbits. Out of the 64 orbits of cones in the chamber complex $\mathcal{G}$, at least 22 orbits (when there are 213 domains) and at most 31 orbits (when there are 337 domains) appear at a time.  

\subsection{Further observations on $A_3$}

With the chamber complex and all the weight polynomials for $A_3$ in hand, we can test whether the bounds on the number of parallel linear factors from Theorems~\ref{thm:factorization} and~\ref{thm:InternalFactorization} are tight in this case. For $k=4$, and a generic dominant weight, we should be getting at least two parallel linear factor in the directions conjugate to $\omega_1$ or $\omega_3$, and three parallel linear factors in the directions conjugate to $\omega_2$, for both the weight polynomials in the boundary regions and the jumps between adjacent regions. 

For generic $\lambda$, two full dimensional domains are adjacent if the corresponding cones in the chamber complex are adjacent. This allows us to test for factorizations on the chamber complex level and thus work with all $\lambda$ at once, so to speak. We have verified that the bounds are met exactly here: when extra linear factors occur, they are not part of the parallel family

For nongeneric $\lambda$, we have to be careful when relating adjacency between regions of the permutahedron and adjacency between cones. This is because, in three dimensions for example, it is possible for two cones to touch along an edge but not along a face, and if we cut them with a hyperplane going through that edge, their two-dimensional intersections become adjacent because they have an edge in common. We get problems on the wall $\lambda_2=\lambda_3$ of the fundamental Weyl chamber, for instance. The weight polynomials in this case all have the form $\gamma(\gamma+1)$, and the jumps between adjacent cones always have the form
\begin{displaymath}
\gamma(\gamma+1) - \gamma'(\gamma'+1) = (\gamma-\gamma')(\gamma+\gamma'+1)\,.
\end{displaymath}
Hence we don't get parallel linear factors at all in this case. For the cases $\lambda_1=\lambda_2$ and $\lambda_3=\lambda_4$, we get 49 or 61 domains generically (see \cite{GLS} for a study of the Duistermaat-Heckman function and its jumps for one of the 49 domain generic cases). The weight polynomials are at most quadratic in this case, and we verify that we get two parallel factors in every jump.

The zero weight does not always appear in the weight space decomposition of $\lambda$, but for $\lambda$ generic, there will be a non-empty domain of the permutahedron that contains the origin (even if it is not a weight). This domain is invariant under the action of the Weyl group and we will call it the \emph{central domain}. 

\begin{figure}[!ht]
\begin{center}
\includegraphics[width=0.6\textwidth]{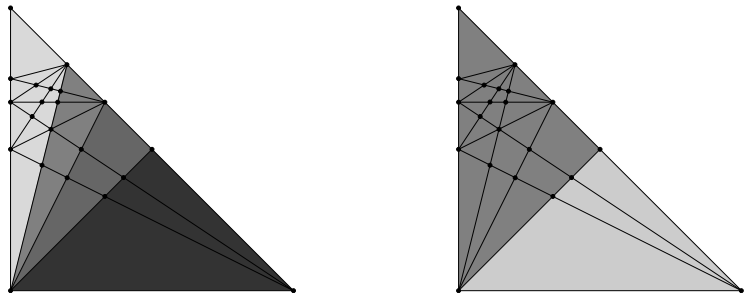}
\caption{Central domains according to the shape of the polytope (left) and the weight polynomials (right).}
\label{fig:CentralDomains}
\end{center}
\end{figure}

We can describe the generic central domains for $A_3$. The diagram on the left in Figure~\ref{fig:CentralDomains} shows the four types of domains for $\lambda_1<-\lambda_4$. The light region corresponds to cubic central domains. In the region next to it, we get truncated cubes: four vertices in tetrahedra position in a cube are truncated, and we get six hexagonal faces and four triangular faces, forming a polytope with 16 vertices. In the remaining two darker regions, we get tetrahedra, in two different orientations (the central domain vanishes on the wall between them). When the hyperplane supporting a wall goes thought the origin, the bounded part giving the wall also contains the origin, and the central domain then vanishes (degenerates to a single point). This occurs when $\lambda=-\lambda^{\mathrm{rev}}$, $\lambda_2=0$ or $\lambda_3=0$.

The behavior of the weight polynomials, as polynomials in both $\lambda$ and $\beta$, is slightly different. For generic $\lambda$, we get only two of them (four if we don't work modulo the Dynkin diagram symmetry). The diagram on the right in Figure~\ref{fig:CentralDomains} shows the separation. The polynomial in the light region does not depend on the $\beta$-coordinates, thus showing that permutahedra for $A_3$ obtained by perturbing around the permutahedra for the fundamental weights $\omega_1$ and $\omega_3$ have tetrahedral central domains over which the multiplicity function is constant (for fixed $\lambda$). Domains like these are called \emph{lacunary} domains in \cite{GLS}. For $\lambda_1<-\lambda_4$, the polynomials are
\begin{eqnarray*}
p^{(\mathrm{light})} & = & \frac{1}{2}(\lambda_2-\lambda_3+1)(\lambda_1-\lambda_2+1)(\lambda_1-\lambda_3+2)\\
p^{(\mathrm{dark})} & = & \frac{1}{2}(\lambda_1-\lambda_2+1)(-\lambda_2^2-2\lambda_3^2+\lambda_3\lambda_4-\lambda_2\lambda_3-\lambda_2\lambda_4+\lambda_2-\lambda_4+2-2h_2(\beta_1\beta_2\beta_3))\,,
\end{eqnarray*}
where $h_2$ is the complete homogeneous symmetric function:
\begin{displaymath}
h_2(\beta_1,\beta_2,\beta_3) = \beta_1^2+\beta_2^2+\beta_3^2+\beta_1\beta_2+\beta_2\beta_3+\beta_1\beta_3\,.
\end{displaymath}

\begin{remark}
Dealing with $A_4$ is more difficult computationally. For example, the permutahedron for the weight $\lambda=\delta$ splits into 15230 regions, and this number is a lower bound on the number of maximal cells of the chamber complex for the weight multiplicities.
\end{remark}

\bigskip

\begin{center}
{\large \textbf{Acknowledgments}}
\end{center}
\medskip

We would like to thank Ron King, Bertram Kostant, Jes\'us De Loera, Richard Stanley, Rekha Thomas, Mich\`ele Vergne and Norman Wildberger for useful discussions and comments.

\medskip

\bibliographystyle{plain}

\bibliography{biblio}

\end{document}